\documentclass[a4paper,12pt]{article}

\pdfoutput=1

\usepackage{amsmath}
\usepackage{amsfonts}
\usepackage{amssymb}
\usepackage[latin9]{inputenc}

\numberwithin{equation}{section}

\usepackage{amsmath}
\usepackage{amsfonts}
\usepackage{amssymb}
\usepackage[latin9]{inputenc}
 \usepackage{amsthm}
\usepackage[all]{xy}

\newtheorem{theorem}{Theorem}[section]
\newtheorem{lemma}{Lemma}[section]
\newtheorem{proposition}{Proposition}[section]
\newtheorem{example}{Example}
\newtheorem{corollary}{Corollary}[section]
\newtheorem{remark}{Remark}[section]

\usepackage{hyperref}

\newcommand{\ga}{\alpha}
\newcommand{\gb}{\beta}

\renewcommand{\ge}{\varepsilon}

\newcommand{\GD}{\Delta}
\newcommand{\SA}{{\cal A}}
\newcommand{\tSA}{\widetilde{\cal A}}
\newcommand{\ot}{\otimes}
\newcommand{\gd}{\delta}

\newcommand{\SC}{{\cal C}}

\renewcommand{\SS}{{\cal S}}

\newcommand{\SB}{{\cal B}}
\newcommand{\gl}{\lambda}

\newcommand{\SI}{{\cal I}}
\newcommand{\SW}{{\cal W}}
\newcommand{\SD}{{\cal D}}

\newcommand{\SV}{{\cal V}}

\newcommand{\SM}{{\cal M}}
\newcommand{\STP}{{\cal M}}
\newcommand{\STPL}{{\cal M}_d}
\newcommand{\STOP}{{\cal OM}}
\newcommand{\STOPL}{{\cal OM}_d}

\newcommand{\RT}{{\rm T}}
\newcommand{\RS}{{\rm S}}
\newcommand{\RH}{{\rm CBH}}
\newcommand{\rest}{\restriction}
\newcommand{\gs}{\sigma}
\renewcommand{\phi}{\varphi}
\renewcommand{\exp}{{\rm exp}}
\newcommand{\GO}{\Omega}
\newcommand{\SL}{{\cal L}}

\renewcommand{\phi}{\varphi}
\renewcommand{\exp}{{\rm exp}}

\newcommand{\id}{{\rm id}}
\newcommand{\bthm}{\begin{theorem}}
\newcommand{\ethm}{\end{theorem}}
\newcommand{\eprop}{\end{proposition}}
\newcommand{\bprop}{\begin{proposition}}
\newcommand{\bexam}{\begin{example}}
\newcommand{\eexam}{\end{example}}
\newcommand{\blem}{\begin{lemma}}
\newcommand{\elem}{\end{lemma}}
\newcommand{\bremark}{\begin{remark}}
\newcommand{\eremark}{\end{remark}}
\newcommand{\bcorol}{\begin{corollary}}
\newcommand{\ecorol}{\end{corollary}}

\newcommand{\bn}{\par \bigskip \par\noindent}

\newcommand{\pn}{\par\noindent}
\newcommand{\balo}{\begin{align*}}
\newcommand{\ealo}{\end{align*}}
\newcommand{\bal}{\begin{align}}
\newcommand{\eal}{\end{align}}
\newcommand{\beq}{\begin{equation}}
\newcommand{\eeq}{\end{equation}}
\newcommand{\beqo}{\begin{equation*}}
\newcommand{\eeqo}{\end{equation*}}

\newcommand{\RA}{{\rm A}}
\newcommand{\RG}{{\rm G}}
\newcommand{\RD}{{\rm D}}

\renewcommand{\mit}{ \, \vert \, }
\newcommand{\SN}{{\cal N}}

\newcommand{\SE}{{\cal E}}
\newcommand{\SF}{{\cal F}}
\newcommand{\RL}{{\rm L}}

\newcommand{\ei}{{\bf 1 }}

\newcommand{\re}{{\rm e}}

\newcommand{\wt}{\widetilde}

\newcommand\mycom[2]{\genfrac{}{}{0pt}{}{#1}{#2}}

\newcommand{\mc}{{\mathbb C}}
\newcommand{\mn}{{\mathbb N}}
\newcommand{\ma}{{\mathbb A}}

\newcommand{\mk}{{\mathbb K}}

\newcommand{\RP}{{\mathbb P}}

\newcommand{\bF}{{\bf F}}

\newcommand{\Sqcup}{\amalg}

\newcommand{\anz}{\#}
\newcommand{\tphi}{\widetilde{\phi}}
\renewcommand{\Sqcup}{\, \underline{\sqcup} \,}

\newcommand{\es}{\emptyset}

\renewcommand{\P}{{\rm P}}

\newcommand{\lan}{\langle}
\newcommand{\ran}{\rangle}
\newcommand{\set}{{\rm set}\,}
\newcommand{\wh}{\widehat}

\newcommand{\boole}{\diamondsuit}

\newcommand{\AlgProb}{{\bf \rm {AlgP}}}
\newcommand{\AlgProbL}{{\bf \rm {AlgP}}_d}
\newcommand{\AlgProbLB}{{\bf \rm { AlgP}}_{d,m}}
\newcommand{\UnitalAlgProbLB}{{\bf \rm { UAlgP}}_{d,m}}
\newcommand{\UnitalAlgProbL}{{\bf \rm { UAlgP}}_d}

\newcommand{\UnitalAlgProb}{\bf \rm { UAlgP}}
\newcommand{\ComUnAlgProb}{\bf \rm { CUAlgP}}

\newcommand{\KC}{{\mathfrak C}}
\newcommand{\KE}{{\rm E}}

\newcommand{\Ten}{(\ot _q) _{q \in \mcd} }
\newcommand{\Boole}{(\boole _{r, s})_{(r, s) \in \mc ^2 \setminus \{(0, 0) \} }}
\newcommand{\Free}{(\circledast _q )_{q \in \mcd}}
\newcommand{\Mon}{(\lhd _q )_{q \in \mcd}}
\newcommand{\Antimon}{(\rhd _q )_{q \in \mcd}}
\newcommand{\mcd}{\dot{\mathbb C }}

\newcommand{\sqcupEi}{\sqcup _{\ei}}

\newcommand{\conv}{\star}

\renewcommand{\i}{\iota}

\newcommand{\mcn}{\mc \langle X_n \rangle }
\newcommand{\hmcn}{\mc \wh{\langle X_n \rangle}}
\renewcommand{\S}{{\rm S}}
\newcommand{\T}{{\rm T}}
\newcommand{\ttheta}{\wt{\theta}}

\newcommand{\End}{{\rm End}}
\newcommand{\Hom}{{\rm Hom}}

\begin{document}

 \title{Non-Commutative Stochastic Independence and Cumulants} 
 \author{ Sarah Manzel$^{\ast}$ \and Michael Sch\"urmann\footnote{The authors acknowledge support by DAAD and MAEDI/MENESR through the PROCOPE programme.}  
 \\
Department of Mathematics and Computer Science  \\
University of Greifswald  \\
Walther-Rathenau-Stra\ss e 47  \\
17487 Greifswald, Germany   \\
}
\date{\today}

\maketitle

\begin{abstract}
In a central lemma we characterize 
\lq generating functions\rq \ of certain functors on the category of algebraic non-commutative probability spaces. Special families of such generating functions correspond to \lq unital, associative universal products\rq \ on this category, which again define a notion of non-commutative stochastic independence.
Using the central lemma, we can prove the existence of cumulants and of \lq cumulant Lie algebras\rq \ for all independences coming from a unital, associative universal product. 
These include the five independences (tensor, free, Boolean, monotone, anti-monotone) appearing in N. Murakis classification, c-free independence of M. Bo\.{z}ejko and R. Speicher, the indented product of T. Hasebe  and the bi-free independence of D. Voiculescu.
We show how the non-commutative independence can be reconstructed from its cumulants and cumulant Lie algebras.
\end{abstract}
    
\section{Preliminaries}
\label{pres}
Vector spaces will be over the field of complex numbers.
Algebras will be complex associative.
An algebra is called unital if it has a unit element.
An algebra with an involution is called a $*$-algebra.
A positive, normalized linear functional on a unital $*$-algebra is called a state.
We say that a linear functional $\phi$ on a $*$-algebra $\SA$  is \emph{strongly positive} if
the linear functional $\tphi$  on $\tSA$ is a state where $\tphi (\ei ) = 1$ and where $\tSA = \mc \ei \oplus \SA$ denotes the unitalization of $\SA$.
\par
Let $\SV$ be a vector space. The tensor algebra over $\SV$ is the vector space
\begin{equation*}
\T (\SV ) = \bigoplus_{n = 1} ^{\infty} \SV ^{\ot n}
\end{equation*}
where $\SV ^{\ot n}$ denotes the $n$-fold tensor product of the vector space $\SV$ with itself. The multiplication of $\T (\SV )$ is given by linear extension of 
\begin{align*}
&(v_1 \ot \ldots \ot v_n ) (w_1 \ot \ldots \ot w_m )  \\
&\quad  = v_1 \ot \ldots \ot v_n \ot w_1 \ot \ldots \ot w_m .
\end{align*}
The algebra $\T (\SV )$ is characterized by the following universal property. For every algebra $\SA$ and linear mapping $R : \SV \to \SA$ there is a unique algebra homomorphism $\T (R) : \T (\SV ) \to \SA$ such that $R = \T (R) \circ \i$ where $\i$ denotes the canonical embedding of $\SV = \SV ^{\ot 1}$ into $\T (\SV )$.
Similarly, the symmetric tensor algebra $\S (\SV )$
\begin{equation*}
\S (\SV ) = \bigoplus_{n = 0} ^{\infty} \SV ^{\ot _s n}
\end{equation*}
 over $\SV$ is defined where now we put $\SV ^{\ot _s 0} = \mc$ and where $\SV ^{\ot _s n}$ denotes the $n$-fold symmetric tensor product of $\SV$ with itself. The commutative unital algebra $\S (\SV )$ is characterized by its universal property, which says that for each \emph{commutative unital} algebra $\SA$ and each linear mapping $R : \SV \to \SA$ there is a unique unital algebra homomorphism $\S (R ) : \S (\SV ) \to \SA$ such that $R = \S (R) \circ \i$.
 Notice that with our notation $\S (\SV )$ is unital whereas $\T (\SV )$ is not.
\par
 Let $I$ be a non-empty (index) set. Denote by $\ma (I)$ the set of all finite sequences $\ge = (\ge _1 , \ldots , \ge _n )$, $n \in \mn$, where $\ge _j \in I$ with $\ge _j \neq \ge _{j + 1}$, $j \in [n - 1]$.
Write $\vert \ge \vert = n$.
We put 
\begin{equation*}
[m] := \{ 1, \ldots , m \}
\end{equation*}
and
$\ma _m := \ma ([m] )$ for $m \in \mn$.
\par
For vector spaces $\SV _i $, $ i \in I$,  we call the vector space
\begin{equation}
\label{freevec}
\bigsqcup_{i \in I} \, \SV _i := \bigoplus_{\ge \in \ma (I)} \SV _{\ge} ,  \ \SV_{\ge} = \SV _{\ge _1} \ot \ldots \ot \SV _{\ge _n}
\end{equation}
the \emph{free product} of the vector spaces $\SV _i $.
For algebras $\SA _i$, $i \in I$, the vector space $\bigsqcup _{i \in I} \, \SA _i$ is turned into an algebra, the \emph{free product} of the algebras $\SA _i$, if we define the multiplication by
\begin{align*}
&(a_1 \ot \ldots \ot a_n ) ( b_1 \ot \ldots \ot b_m )  \\
&\hspace{0,5cm} = \left\{ \begin{array}{ll} a_1 \ot \ldots \ot a_{n - 1}\ot (a_n b_1 ) \ot b_2 \ot \ldots \ot b_m &\mbox{if } \ge _n = \rho _1   \\
a_1 \ot \ldots \ot a_n  \ot  b_1 \ot \ldots \ot b_m  &\mbox{if } \ge _n \neq \rho _1
\end{array}
\right.
\end{align*}
for $\ge , \rho \in \ma (I)$, $\ge = (\ge _1 , \ldots , \ge _n )$, $\rho = (\rho _1 , \ldots , \rho _m )$, $a_1 \ot \ldots \ot a_n \in \SA _{\ge}$, $b_1 \ot \ldots \ot b_m \in \SA _{\rho}$.
For an algebra $\SA$ we put
\begin{equation}
\label{freeprod}
\SA ^{\sqcup I} :=
\bigsqcup_{i \in I} \SA _i , \ \SA _i \ \mbox{a copy of} \ \SA .
\end{equation}
Let $\SA$ be an algebra. We say that the
$(\SA ^{(1)} , \ldots , \SA ^{(m)} )$, $m \in \mn$, is \emph{free} in $\SA$ if $\SA ^{(1)} , \ldots , \SA ^{(m)}$ are subalgebras of $\SA$ and if the algebra homorphism from $\SA ^{(1)} \sqcup \ldots \sqcup \SA^{(m)}$ to $\SA$ which is given by
\begin{equation*}
\SA _{\ge} \ni a_ 1 \ot \ldots \ot a_n \mapsto a_1 \ldots a_n \in \SA
\end{equation*}
is injective. If this mapping is a bijection, we say that
$ (\SA ^{(1)} , \ldots , \SA ^{(m)})$ \emph{freely generates} $\SA$.
In this case, we also write $\SA = \SA ^{(1)} \sqcup \ldots \sqcup \SA^{(m)}$.
\par
For a non-empty set $X$ we denote by $\lan X \ran$ the \emph{semi-group} freely generated by $X$ and by $[X]$ the \emph{commutative monoid} freely generated by $X$.
Notice that, with our notation, $\langle X \rangle$ does not include the empty word whereas $[X]$ does.
The elements of $\lan X \ran$ are formal words $x_1 \ldots x_k$, $k\in \mn$, $x_i \in X$, which are called monomials.
The monomials of $[X]$ will sometimes be written $x_1 \cdot \ldots \cdot x_k$, $k \in \mn _0$ with $\mn _0 = \mn \cup \{ 0\}$.
We will deal with monomials of $[\langle X \rangle ]$, which are of the form
$M_1 \cdot \ldots \cdot M_l$, $l \in \mn _0$,  with $M_i = x_{i,1} \ldots x_{i, k_i}$, $i \in [l]$, $k_i \in \mn$.
We put
\begin{equation*}
\set ( x_1 \ldots x_k ) = \{x_1 , \ldots , x_k \}
\end{equation*}
 and 
\begin{equation*}
\vert x_1 \ldots x_k \vert = k .
\end{equation*}
For a set $X$ denote by $\mc X$ the vector space with basis $X$ and by $\mc \lan X \ran $ the \lq semi-group algebra\rq \ of $\lan X \ran$, i.e. the vector space with basis $\lan X \ran$ and algebra multiplication given by the bilinear extension of the multiplication of $\lan X \ran$. Then $ \T (\mc X ) \cong \mc \langle X \rangle$ in a natural way.
Similarly, the (commutative, unital) polynomial algebra $\mc [X]$ will be considered. 
We have $\S (\mc X ) \cong \mc [X]$.
We will also use the natural identifications 
\begin{equation}
\RS (\SV _1 \oplus \SV _2 ) = \RS (\SV _1 ) \ot \RS (\SV _2 )
\end{equation}
and
\begin{equation}
\label{freesum}
\RT (\SV _1 \oplus \SV _2 ) = \RT (\SV _1 ) \sqcup \RT (\SV _2 ) .
\end{equation}
\par
Let $\SV ^*$ denote the algebraic dual space of a vector space $\SV$.
For $\phi \in \SV ^*$ we have the scalar-valued algebra homomorphism $\S (\phi ) : \S (\SV ) \to \mc$, and we say that a mapping $F : \SV ^* \to \mc$ is a \emph{polynomial function} over $\SV$ if there is an element $P \in \S (\SV )$ such that $F (\phi ) = \S (\phi ) (P)$ for all $\phi \in \SV ^*$. In fact, the polynomial $P = P _F$ is uniquely determined by its polynomial function $F$.
If $\{ v_i \mit i \in I \}$ is a vector space basis of $\SV$, then  we can identify $\S (\SV )$ and $\mc [\{x_i \mit i \in I \}]$ and we have $\S (\phi ) (P ) = P ((\phi (x_i )_{i \in I })$, that is $\S (\phi ) (P )$ is the value of the polynomial function given by $P$ at the \lq point\rq \ $(\phi (x_i ))_{i \in I}$.
\par
Let $\SV$ and $\SW$ be vector spaces and $d \in \mn$.
Frequently, we use the identifications
\begin{equation*}
\Hom_{\mc} (\SV ^d , \SW ) = \bigl( \Hom _{\mc} (\SV , \SW ) \bigr) ^d = \Hom _{\mc} (\SV , \SW ^d ).
\end{equation*}
For example, $T v = (T^{(1)} v , \ldots , T^{(d)} v ) \in \SW ^d$ or 
$T (v^{(1)} , \ldots , v^{(d)} ) = T^{1)} v^{(1)} + \ldots + T^{(d)} v^{(d)} \in \SW$ for 
\begin{equation*}
(T^{(1)} , \ldots , T^{(d)} ) \in \bigl( \Hom_{\mc} (\SV , \SW )\bigr) ^d, v, v^{(1)}, \ldots , v^{(d)}  \in \SV .
\end{equation*}
\bn
Put $X_n  := \{ x_ 1, \ldots , x_ n \}$, the set formed by $n \in \mn$ \lq indeterminates\rq \ $x_1, \ldots , x_n $.
Given $n \in \mn$ elements $a_1, \ldots , a_n$ of an algebra $\SA$ we let
\begin{equation*}
j(a_1 , \ldots , a_n ) : \mc \lan X_n \ran \to \SA
\end{equation*}
be the algebra homomorphism defined by putting
\begin{equation*}
j(a_i ) = x_i , \ i \in [n] .
\end{equation*}
Considering $\mc \lan X \ran = \RT (\mc X )$ we have that  
$$
j (a_1 , \ldots , a_n )= \RT (R)
$$ 
with $R : \mc X \to \SA$ linear, $R(x_i ) = a_i$.
Put ($d \in \mn$) 
\begin{equation*}
\SV ^{*d} : = (\SV ^* ) ^d = \SV ^* \oplus \ldots \oplus \SV ^* = (\SV \oplus \ldots \oplus \SV ) ^* = (\SV ^d )^*
\end{equation*}
and write $\phi \circ j$ for $(\phi ^{(1)} \circ j , \ldots , \phi ^{(d)} \circ j )
\in \SW ^{*d} $ if 
\begin{equation*}
\phi = (\phi ^{(1)} , \ldots , \phi ^{(d)} ) \in \SV ^{*d}
\end{equation*}
 where $\SV$, $\SW$ are vector spaces and $j : \SW \to \SV$ is linear.
\par
The polynomial algebra $\mc [\langle X \rangle ]$ in the indeterminates $M$, $M \in \langle X \rangle$, is an $\mn _0$-graded algebra if we put the degree $\vert M_1 \cdot \ldots \cdot M_l \vert$ of a monomial $M_1 \cdot \ldots \cdot M_l$ in $[\langle X \rangle ]$ equal to $\vert M_1 \vert + \ldots + \vert M_l \vert$.
For a finite set $X$, $\# X = n$, denote by $\STP (X)$ the set of all monomials in $\mc [ \langle X \rangle ]$ of the form $M_1 \cdot \ldots \cdot M_l$ with
\begin{gather*}
\vert M_1 \cdot \ldots \cdot M_l \vert = n  \\
(\set M_r ) \cap(\set M_s ) = \emptyset \ \mbox{for} \ r \neq s  \\
\# \set M_r = \vert M_r \vert , r \in [l] .
\end{gather*}
The set $\STP (X)$ is formed by all elements of $[\langle X \rangle ]$ with the property that each indeterminate of $X$ appears exactly once. 
For example, for $X = X_4$, we have that 
\begin{equation*}
(x_3 x_1 x_4) \cdot x_2 , x_4 \cdot x_1 \cdot x_3 \cdot  
x_2 , \ x_2 x_1 x_4 x_3  \mbox{ and }  x_1 \cdot x_2 \cdot (x_4 x_3 )
\end{equation*}
are in $\SM (X)$
whereas 
\begin{equation*}
(x_1 x_2 ) \cdot x_4 \mbox{ or } (x_1 x_2 x_4) \cdot (x_1 x_3 )
\end{equation*}
 are not.
Denote by $\STPL (X)$ the set 
\begin{align}
\label{labelled}
&\{ (M _1 , k _1 ) \cdot \ldots \cdot (M_l , k_l )   \mit M _1 \cdot \ldots \cdot M_l \in \STP (X) ,  \\
\nonumber
&\hspace{6cm}  k_1 , \ldots , k_l  \in [d] \}
\end{align}
 of \lq labelled elements\rq \  of $\STP (X)$, which is a subset of $[\langle X \rangle \times [d] ]$, and put
\begin{equation*}
\STP (n) = \STP (X_n ), \ \STPL (n) = \STPL (X_n ) .
\end{equation*}

\section{Introduction}
Let $(\KC , \odot, \KE )$ be a tensor category; see, for example, \cite{McL} or \cite{AgMa}.
We will always assume that the unit object $\KE$ is an initial object of the category $\KC$. We write $\ei _A : \KE \to A$ for the unique morphism from $E$ to an object $A$ of $\KC$. Denote by 
\begin{gather*}
\ga _{A, B, C} : A \odot (B \odot C) \to (A \odot B ) \odot C  \\
l _A : \KE \odot A \to A , \  r _A : A \odot \KE \to A
\end{gather*}
the associativity and unit constraints of the given tensor category, which satisfy the pentagon and triangle diagrams.
In  this situation the morphisms $\iota _1 : A_1 \to A_1 \odot A_2 $ and $\iota _2 : A_2 \to A_1 \odot A_2 $ with
\begin{align*}
\iota _1 &=  (\id _{A_1} \odot \ei _{A _2})\circ r_{A_1}^{-1}  \\
\iota _2 &=  (\ei _{A _1} \odot \id _{A_2} )\circ l_{A_2}^{-1}
\end{align*}
are called the inclusion maps.
Following U. Franz \cite{Fra06}, M. Gerhold \cite{Ger1} and S. Lachs \cite{Lachs1}, two morphisms $j_1 : A_1 \to A$, $j_2 : A_2 \to A$ with the same target are called \emph{independent} if there is a unique morphism $h : A_1 \odot A_2 \to A$ such that $j_1 = h \circ \iota _1 $ and $j_2 = h \circ \iota _2$.
For example, the tensor product $\odot$ is a coproduct for the category $\KC$ if and only if \emph{all} pairs of morphisms with the same target are independent.
As another elementary example consider the category of vector spaces over a field $\mk$ with injective linear mappings as morphisms and the tensor product given by the direct sum. Then two morphisms with the same target are independent if and only if there images are linearly independent; see \cite{Ger1}.
\par
Now consider the category of duals of probability spaces. The objects of this category are pairs $(\SN , \phi )$ with $\SN$ a commutative von Neumann algebra and $\phi$ a normal state on $\SN$. The morphisms are von Neumann algebra homomorphisms which respect the states. A random variable $X: \GO \to E$ over a probability space $(\GO , \SF , \RP )$ with values in a measurable space $(E, \SE )$ gives rise to the morphism 
\begin{equation*}
j_X : \RL ^{\infty } (E , \SE , \RP _X ) \to \RL ^{\infty } (\GO , \SF , \RP ); \
j_X (f ) = f \circ X
\end{equation*}
where $\RP _X$ denotes the law of $X$.
The von Neumann algebra tensor product and the tensor product of states turn this category into a tensor category with unit object $\mc$, which is an initial object.
Two random variables $X _1 : \GO \to E_1$, $X_2 : \GO \to E_2$ are stochastically independent in  the usual sense of classical probability theory if and only if the morphisms $j_{X_1}$, $j_{X_2}$ are independent in the above sense.
\par
In this paper we will study non-commutative notions of independence, and, in particular, the construction of cumulants associated with such an independence. A natural approach to non-commuta\-ti\-vity would be to consider the category of von Neumann algebras and states as above without the assumption of commutativity. Then one could ask for possible tensor products on this category. The classification of tensor products in this analytic setting seems to be an open problem. There are classification results for a more algebraic approach. 
Consider the category $\AlgProb$ with pairs $(\SA , \phi )$, $\SA$ an algebra, $\phi$ a linear functional on $\SA$, as objects, and morphisms given by algebra homomorphisms which respect the linear functionals. The initial object is given by the algebra $\{ 0 \}$ and the only linear functional $0 : \{ 0\} \to \mc$ on $\{ 0 \}$. 
Let
$(\AlgProb, \odot ' , \{ 0\} )$ be a tensor category and assume that
\begin{equation*}
(\SA _1 , \phi _1 ) \odot ' (\SA _2 , \phi _2 ) = (\SA _1 \odot ' \SA _2 , \phi _1 \odot ' \phi _2 ) .
\end{equation*}
We form a new tensor category $(\AlgProb , \odot , \{ 0 \} )$ by putting
\begin{gather*}
\SA _1 \odot \SA _2 = \SA _1 \sqcup \SA _2 \\
\phi _1 \odot \phi _2 = (\phi _1 \odot ' \phi _2 ) \circ (\iota _1 \odot ' \iota _2 )
\end{gather*}
where here $\iota _1$, $\iota _2$ are the embeddings 
 of the algebras $\SA _1$, $\SA _2$ into $\SA _1 \odot ' \SA _2$ given 
 by $\odot '$ as above. 
The free product is a coproduct in the category of algebras. The tensor product $\odot$ is a reduction of the tensor product $\odot '$ in the sense of \cite{Fra06}. If two morphisms are independent with respect to $\odot '$, they are also independent with respect  to $\odot$.
This means that, for many questions in non-commutative probability,  we may restrict ourselves to tensor products of the form $(\SA _1 \sqcup \SA _2 , \phi _1 \odot \phi _2 )$ where the product $\odot$ of linear functionals satisfies
\begin{gather}
\label{A1}
(\phi _1 \odot \phi _2 ) \circ \i _1 = \phi _1 ; \ (\phi _1 \odot \phi _2 ) \circ \i _2 = \phi _2  \\
\label{A2}
(\phi _1 \odot \phi _2 ) \odot \phi _3 = \phi _1 \odot (\phi _2 \odot \phi _3 ) \\
\label{A3}
(\phi _1 \odot \phi _2 ) \circ ( j_1 \,  \underline{\sqcup} \, j_2 ) = (\phi _1 \circ j_1 ) \odot (\phi _2 \circ j _2 )
\end{gather}
for objects $(\SA _1 , \phi _1 )$, $(\SA _2 , \phi _2 )$ and $(\SA _3 , \phi _3 )$ in $\AlgProb$, algebras $\SB _1$, $\SB _2$, and algebra homomorphisms $j_1 : \SB _1 \to \SA _1$, $j_2 : \SB _2 \to \SA _2$.
Here we put
\begin{equation*}
j_1 \, \underline{\sqcup } \, j_2 = (\iota _1 \circ j_1 ) \sqcup (\iota _2 \circ j_2 ) .
\end{equation*}
Products in $\AlgProb$ which satisfy (\ref{A1})-(\ref{A3}) have been completely classified; see \cite{Spe97,BGS02,Mur02,Mur03,GeLa}.  There are five families of such products, the tensor \cite{HudPar}, free \cite{Voi85,Sp} , Boolean \cite{vW1}, monotone and anti-monotone \cite{Mur97,Mur01,Lu} families.
\par
In this paper we will consider a more general situation. The category $\AlgProbLB$, $d, m \in \mn$, has as objects triplets $\bigl( \SA , (\SA ^{(1)} , \ldots , \SA ^{(m)} ), \phi \bigr)$ where now $\SA$
comes with an $m$-tuple $(\SA ^{(1)} , \ldots , \SA ^{(m)} )$ of subalgebras which freely generates $\SA$,
 and $\phi$ 
 is a $d$-tuple of linear functionals with each
 component 
a linear functional on  $\SA$.
A morphism of $\AlgProbLB$ is given by an $m$-tuple of algebra homomorphisms adapted to the free decomposition of the algebras into subalgebras and which  respect the components of $\phi$.
\par
We consider tensor products  on $\AlgProbLB$ given by products $\odot$ with (\ref{A1})-(\ref{A3}) which now have to be interpreted component-wise.
We  call such a unital, associative universal product a \emph{u.a.u.-product}.
Now further examples of independences, like bi-freeness \cite{Voi123} or c-freeness \cite{BSp,Ansh} and the examples of \cite{Has}, are included. The classification of u.a.u.-products in $\AlgProbLB$ is an open problem.
\par
We are interested in convolution products, moments and cumulants in the non-commutative algebraic setting. To define these notions some additional structure is needed. Classically, convolution of measures on (topological) groups can be defined by a combination of the tensor product of two measures and the group multiplication.
Similarly, the convolution of linear functionals on a Hopf algebra or bialgebra is defined.
Since we wish to study different notions of tensor products in a non-commutative setting, we need \emph{dual groups} or semi-groups \emph{in the sense of D. Voiculescu} \cite{Voi87}, which better suit non-commutativity than Hopf algebras or bialgebras.
A dual semi-group in this paper means a pair $(\SD , \GD )$ with $\SD$ an algebra and $\GD : \SD \to \SD \sqcup \SD$ an algebra homomorphism such that
\begin{align}
(\GD \Sqcup \id ) \circ \GD &= (\id \Sqcup \GD ) \circ \GD \\
(0 \Sqcup \id ) \circ \GD &= \i _1 ; \
(\id \Sqcup 0 ) \circ \GD = \i _2 .
\end{align}
\par
For the classification of independences a lemma from \cite{BGS05} is crucial.
A generalized version of this lemma and a more detailed proof can be found in Section \ref{basiclemma} of this paper.
\par
We explain the statement of our central lemma.
The universal property of the polynomial algebra $\RT (\mc X_ n ) = \mc \langle X_n \rangle$ can be formulated as follows.
Given an arbitrary algebra $\SA$ and $n$ elements $a_1 , \ldots , a_n$ in $\SA$ there is a unique algebra homomorphism $j(a_1 , \ldots , a_n ) : \mcn \to \SA$ which maps $x_i$ to $a_i$, $i \in [n]$. 
A mapping $F : \mcn ^{*d} \to \mc$ from the space $\mcn ^{*d}$ of $d$-tuples of linear functionals on $\mcn$ to the complex numbers is called \emph{generating} if the maps
\begin{equation*}
(a_1 , \ldots , a_n ) \to F(\phi \circ j(a_1 , \ldots , a_n ))
\end{equation*} 
are $n$-linear for all algebras $\SA$ and all $d$-tuples $\phi$ of linear functionals  on $\SA$.
In the special case $d = 1$ the central lemma says that generating maps are always polynomial functions given by a polynomial $P_F$ of degree $n$ in the polynomial algebra $\mc [\langle X_n \rangle ]$ in \emph{commuting} indeterminates $M$ where $M$ runs through all monomials of $\mcn$.
Moreover $P_F$ is of very special type. It lies in the linear span of elements of the form $M_1 \cdot \ldots \cdot M_l$ with monomials $M_k$ in $\mcn$ such that each indeterminant $x_1 , \ldots , x_n$ appears exactly once in $M_1 , \ldots , M_l$, i.e. $P$ is in the linear span of $\SM (X_n )$ of Section \ref{pres}.
$F$ is obtained from $P_F$ by the equation $F (\phi ) = \P_F \bigl( (\phi (M))_M \bigr)$ that is by the value of the  \emph{polynomial function} $P_F$ at the points $\phi (M)$.
For general $d \in \mn$ the monomials $M$ get a label $j \in [d]$ and the $d$-tuple $\phi$ puts the value $\phi ^{(j)} (M)$ to the labelled monomial $(M, j)$.
\par
The same lemma is used to define the \lq Lachs functor\rq \ (which is an example of a colax monoidal functor; see e.g. \cite{AgMa}), the case $d = m = 1$ was treated by S. Lachs in \cite{Lachs1}. Fix a u.a.u.-product $\odot$ on $\AlgProbLB$. The convolution product of two linear functionals $\phi _1$, $\phi _2$ on a dual semi-group $\SD$ is defined by 
\begin{equation}
\phi _1 \conv \phi _2 = (\phi _1 \odot \phi _2 ) \circ \GD .
\end{equation}
The  product $\odot$ gives rise to a cotensor (or colax monoidal) functor $\SL _{\odot}$ from the tensor category $(\AlgProbLB , \odot , \{ 0 \} )$ to the tensor category of commutative algebraic probability spaces  with the usual tensor product of algebras and of normalized linear functionals. It associates with each dual semi-group $\SD$ a bialgebra structure on the symmetric tensor algebra $\SL _{\odot} (\SD ) = \S ((\SD ^{\sqcup m})^d)$ over the vector space $(\SD ^{\sqcup m})^d$ such that 
\begin{equation}
\S (\phi _1 \conv \phi _2 ) = \S (\phi _1 ) \conv \S (\phi _2 )
\end{equation}
where the second convolution is given by the coproduct $\SL _{\odot} (\GD )$ of the bialgebra $(\SL _{\odot} (\SD ), \SL _{\odot} (\GD ), \S (0))$ and $\S (\phi )$ denotes the algebra homomorphism from $\SL _{\odot} (\SD )$ to the complex numbers given by $\phi \in (\SD ^{\sqcup m})^{*d}$.
If $\SD$ is a dual group, the vector space $(\SD ^{\sqcup m})^{*d}$ is a group with respect to convolution; cf. \cite{FrMcK} for the case $d = m = 1$.
\par
Now consider the following example of a dual group which is important for the treatment of cumulants associated with a given product $\odot$ on $\AlgProbLB$. The algebra structure of this dual group is just that of a tensor algebra $\T (\SV )$ over a given vector space $\SV$. 
As before it is sometimes convenient to look at $\T (\SV )$ as the polynomial algebra $\mc \langle x_i ; i \in I \rangle$ in indeterminates $x_i$ with $I$ a set of the cardinality of a basis of $\SV$.
The dual space of $\mc \langle x_i ; i \in I \rangle$ can be identified with the algebra $\mc \langle\langle x_i ; i \in I \rangle\rangle$ of formal power series in the $x_i$.
The comultiplication of $\RT (\SV )$ is given by $\GD x_i = \i _1 x_i + \i _2 x_i$
where $\i _1 , \i _2$ denote the natural embeddings of $\SV$ into $\T (\SV )\sqcup \RT (\SV ) = \RT (\SV \oplus \SV )$; see \ref{freesum}.
\par
The fact that the polynomial in the central lemma is of the special type described above means that the bialgebra $\SL _{\odot} (\SD )$ in the case $\SD = \RT (\SV )$ is $\mn _0$-graded with a 1-dimensional $0$-component. For such $\mn _0$-graded bialgebras (which are automatically Hopf algebras; see e.g. \cite{Man})  linear functionals which vanish at $\ei$ are nilpotent with respect to convolution if restricted to the sub-coalgebra spanned by homogeneous elements of a degree which is bounded by a fixed natural number. 
Thus we may define, pointwise, the linear functional $\ln _{\conv} \Phi$ 
as a logarithm series with respect to convolution for all linear functionals $\Phi$ on $\SL _{\odot} \bigl( \RT (\SV ) \bigr)$ with $\Phi (\ei ) = 1$.
It can be shown that, 
for $\phi \in (\SD ^{\sqcup m})^{*d}$,  $\SD = \RT (\SV )$,
the convolution logarithm $\ln _{\conv} \S (\phi )$ is an $\S (0)$-derivation.
 As such it is determined by its restriction $\ln _{\odot} (\phi )$ to $(\SD ^{\sqcup m})^{d}$.
We call $\ln _{\odot} (\phi ) \in (\SD ^{\sqcup m})^{*d}$ the \emph{cumulant functional} of $\phi$.
In fact, $d = m = 1$ gives the known cumulants (see e.g. \cite{AHLV}) if we choose for $\odot$ one of the Muraki five \cite{Mur02,Mur02,Lachs1}, and the cumulants in the bi-free  case of Voiculescu \cite{Voi123} for $d = 1$, $m = 2$.
Morever, $m = 1$, $d \geq 2$ covers c-freeness of Bo\.{z}ejko and Speicher \cite{BSp,BLSp} and the examples of Hasebe \cite{Has}.
In the symmetric case, e.g. for tensor, free, Boolean and bi-free independence, we have 
\begin{equation}
\label{Add}
\ln _{\odot} (   \phi _1 \conv \phi _2 ) = \ln _{\odot} (\phi _1 ) + \ln_{\odot} (\phi _2 )
\end{equation}
 whereas, for example, in the monotonic or anti-monotonic case this formula has to be replaced by a Campbell-Baker-Hausdorff expansion; cf. \cite{FrMcK} and Theorem \ref{CBH2} of this paper. The uniqueness of cumulants which satisfy natural conditions is shown in \cite{Has2} in the monotonic case.
\par
The Campbell-Baker-Hausdorff expansion follows easily from the calculus of formal power series. The linear functional $\exp _{\conv}\Psi$ can be defined as an exponential convolution series 
 for all linear functionals $\Psi$ on  $\SL _{\odot} \bigl( \RT (\SV ) \bigr)$ and, if $\Psi$ vanishes at $\ei$, we have that this series is a finite sum at each point.
Moreover, it follows that
$\ln _{\conv}$ and $\exp _{\conv}$ are inverses of each other.
For an $\S (0)$-derivation $D = D (\psi)$, $\psi$ the restriction of $D$ to $\bigl( \RT (\SV )^{\sqcup m} \bigr) ^d = \bigl( \RT (\SV ^{m} ) \bigr) ^d$, we have that $\exp _{\conv} D$ is a unital algebra homomorphism and as such of the form $\S (\phi )$ with $\phi \in \bigl( \RT (\SV^{ m}) \bigr) ^{*d}$.
Put $\exp _{\odot} (\psi ) = \phi$. Then $\ln _{\odot}$ and $\exp _{\odot}$ are inverses of each other.
Moreover,
$\ln_{\odot} (\phi _1 \conv \phi _2 )$ is the restriction to $(\SD ^{\sqcup m})^{d}$ of the Campbell-Baker-Hausdorff expansion
\begin{align*}
&\ln _{\conv} \S (\phi _1 ) + \ln _{\conv } \S (\phi _2 ) +  \frac12 [\ln _{\conv} \S (\phi _1 ), \ln _{\conv} \S (\phi _2 )]   _{\conv} \\
& \hspace{3cm} + \frac{1}{12} [\ln _{\conv} \S (\phi _1 ), 
[\ln _{\conv} \S (\phi _1 ), \ln _{\conv} \S (\phi _2 ) ]_{\conv} ] _{\conv} - ...
\end{align*}
which again is a finite sum at each point.
 The convolution commutator of two $\S (0)$-derivations is again an $\S (0)$-derivation.
Therefore, the vector space $(\SD ^{\sqcup m})^{*d}$, $\SD = \T (\SV )$, becomes a Lie algebra with the Lie bracket defined by
\begin{equation}
[ \psi _1 , \psi _2 ] _{\odot} := [ D(\psi _1 ) , D(\psi _2 ) ]_{\conv} \rest (\SD ^{\sqcup m})^{d}
\end{equation}
which we call the \emph{cumulant Lie algebra} for $\SV$. 
Call the cumulants for $\SV = \mc ^n$, $n \in \mn$, i.e. in the case   $\SD = \mc \langle x_1 , \ldots , x_n \rangle$, the \emph{$n$-th order cumulants} of $\odot$ and the associated Lie algebra on $(\SD ^{\sqcup m})^{*d} = \mc \langle \langle x_1 , \ldots , x_{mn} \rangle \rangle ^d$ the \emph{$n$-th order cumulant Lie algebra} of $\odot$.
We show that a u.a.u.-product is determined by the following data
\pn
\begin{itemize}
\item
all $n$-th order cumulants, $n \in \mn$
\item
all $n$-th order cumulant Lie algebras,  $n \in \mn$ .
\end{itemize}
The system of cumulants and cumulant Lie algebras is the infinitesimal form of the underlying u.a.u.-product.
It will be a subject of further research to find general conditions for an \lq infinitesimal system\rq \ to lead to a u.a.u.-product by \lq exponentiating\rq .
\par
We give a brief description of the content of the paper.
Section 3 contains a generalization of Lemma 1 of \cite{BGS02} to the general $d,m$-case. In Sections 4 and 5 we apply this central lemma to u.a.u.-products to find a coefficient formula which could be a first step to the classification of u.a.u.-products. Moreover, we prove a generalization of a result of N. Muraki \cite{Mur4} which says that for \emph{positive} products all \lq wrong ordered\rq \ coefficients vanish.
Section 6 is about the Lachs functor. Section 7 presents some results on $\mn _0$-graded Hopf algebras which are used to define cumulants in the general $d,m$-case in Section 8.

\section{Universal products}
Consider the category of algebraic probability spaces $\AlgProb$ with pairs $(\SA , \phi )$ as objects and with morphisms $j: (\SB , \psi ) \to (\SA , \phi )$, $j : \SB \to \SA$ an algebra homomorphism, $\phi \circ j = \psi$.
More generally, denote by $\AlgProbL$ the category formed by pairs $(\SA  , \phi)$ with $\phi \in \SA ^{*d}$ and morphisms given by algebra homomorphisms which respect each of the components of $\phi$ and $\psi$, i.e. $\phi ^{(l)} = \psi ^ {(l)} \circ j$ for $l \in [d]$. We will write $\phi \circ j$ for $(\phi ^{(1)} \circ j, \ldots , \phi ^{(d)} \circ j )$ so that the morphism condition reads again $\phi \circ j = \psi$. 
Similarly, we consider the category $\UnitalAlgProbL$ of pairs $(\SA , \phi )$ with $\SA$ a \emph{unital} algebra and $\phi = (\phi ^{(1)} , \ldots, \phi ^{(d)} )$ a $d$-tuple of \emph{normalized} linear functionals on $\SA$, i.e. $\phi ^{(l)} (\ei ) = 1$.
Of course, the morphisms are now the \emph{unital} algebra homomorphisms respecting the $\phi ^{(l)}$. 
\par
A  universal product in $\AlgProbL$ is a bifunctor of the special form
\begin{align}
\label{bi-functor}
&\bigl( (\SA _1 , \phi _1 ), (\SA _2 , \phi _2  ) \bigr)  \\
\nonumber
 & \hspace{1cm} \mapsto 
(\SA _1 , \phi _1 ) \odot (\SA _2 , \phi _2 ) = (\SA _1 \sqcup \SA _2 , \phi _1 \odot \phi _2 )  \\
\nonumber
& (j_1 , j_2 ) \mapsto j_1 \Sqcup j_2 .
\end{align}
The bifunctor property implies that (\ref{A3})
holds. A universal product is called degenerate if
\begin{equation*}
(\phi _1 \odot \phi _2 )(a_1 \ldots a_m ) = 0
\end{equation*}
for all objects $(\SA _1 , \phi _1 ), (\SA _2 , \phi _2 )$ in $\AlgProbL$, $\ge \in \ma _2$, $m = \vert \ge \vert > 1$, $a_1 \ldots a_m \in \SA _{\ge}$.
We always assume that our universal products are non-degenerate.
$\odot$ is a u.a.u.-product if it satisfies (\ref{A1}), (\ref{A2}) and (\ref{A3}).
We say that $\odot$ is {\it symmetric} if
\begin{equation}
\label{A4}
\phi _1 \odot \phi _2 = \phi _2 \odot \phi _1 .
\end{equation}
Universal products in the category $\UnitalAlgProbL$ are defined in an analogous way. 
In this case the role of the free product $\SA _1 \sqcup \SA _2$ of two algebras is played by the coproduct in the category of unital algebras, which is given by $\SA _1 \sqcupEi \SA _2 : = (\SA _1 \sqcup \SA _2 )/ \SI$ where $\SI$ is the ideal generated by $\ei _1 - \ei _2$ with $\ei _1, \ei _2$ the unit elements of $\SA _1, \SA _2$.
We mention that for the unitalization $\wt{\SA} = \mc \ei \oplus \SA$ of an algebra
$\SA$ we have in a natural way
\begin{equation*}
\wt{\SA _1} \sqcupEi \wt {\SA _2} \cong \wt{\SA _1 \sqcup \SA _2 } .
\end{equation*}
We will also consider universal products $\odot _{i \in I} \, \phi _i$ of a \emph{family} 
\begin{equation*}
(\SA _i , \phi _i )_{i \in I}
\end{equation*}
 of objects in $\AlgProbL$ or in $\UnitalAlgProbL$. Condition (\ref{A3}) is replaced by
\begin{equation*}
\bigl( \odot _{i \in I} \, \phi _i \bigr) \circ \bigl( \underline{\sqcup}_{i \in I} \, j_i \bigr) = \odot _{i \in I} (\phi _i \circ j_i ) .
\end{equation*}
Now we generalize from one to $m\in\mn$ \emph{algebras} by considering the category $\AlgProbLB$ with objects $\bigl( \SA, (\SA ^{(1)} , \ldots , \SA ^{(m)} ), \phi \bigr)$ where 
\begin{gather*}
(\SA ^{(1)} , \ldots , \SA ^{(m)} ) \mbox{ freely generates } \SA ,\\
\phi = (\phi ^{(1)} , \ldots , \phi ^{(d)} ) ,\
\phi ^{(1)} , \ldots , \phi ^{(d)} \in \SA ^{*},
\end{gather*}
and morphisms given by algebra homomorphisms $j : \SB \to \SA$ such that $j(\SB ^{(k)}) \subset \SA ^{(k)}$, $k \in [m]$,
and  $\psi = \phi \circ j$.
Define the category $\UnitalAlgProbLB$ in an analogous way.
\par
A universal product in the category $\AlgProbLB$ is a bifunctor of the form
(\ref{bi-functor}) where now $\SA _1 \sqcup \SA _2 $ stands for
\begin{equation*}
\SA _1 \sqcup \SA _2 , \bigl( (\SA _1 ^{(1)} \sqcup \SA _2 ^{(1)} ), \ldots , (\SA _1 ^{(m)} \sqcup \SA _2 ^{(m)} ) \bigr),
\end{equation*}
which means again (\ref{A3}).
As before, we say that $\odot$ is unital and associative  if (\ref{A1}) and (\ref{A2}) hold.
\bexam
\label{five}
{\rm
In the case $d = m = 1$,  by a result of S. Lachs \cite{Lachs1}, there are exactly five families $\Ten$, $\Boole$, $\Free$, $\Mon$, $\Antimon$ of u.a.u.-products. 
The tensor family $\Ten$ can be obtained from the parameter $q = 1$ case by the normalization equation
\begin{equation} 
\label{q-norm}
\phi _1 \otimes _q \phi _2 = q ^{-1} \bigl( (q \, \phi _1 ) \otimes _1 (q \, \phi _2 ) \bigr) .
\end{equation}
Moreover,
\begin{equation*}
\phi _1 \ot _1 \phi _2 = (\phi _1 \ot \phi _2 ) \circ \gs
\end{equation*}
with 
\begin{equation*}
\gs : \SA _1 \sqcup \SA _2 \to \SA _1 \oplus \SA _2 \oplus (\SA _1 \ot \SA _2 )
\end{equation*}
 the canonical algebra homomorphism from the free product to the tensor product of the algebras $\SA _1$ and $\SA _2$.
Similarly, the free family $\Free$ and the monotonic $\Mon$ and anti-monotonic family $\Antimon$ can be obtained from the $q = 1$ case by a normalization equation like in (\ref{q-norm}). 
For $a_1 \ot \ldots \ot a_n \in \SA _{\ge}$, $\ge = (\ge _1 , \ldots , \ge _n ) \in \ma _2$, the free $q = 1$ case  is given by the recursion formula 
\begin{align*}
&(\phi _1 \circledast \phi _2 )(a_1 \ot \ldots \ot a_n )   \\
& \ = \sum _{I \subsetneq [n]} (-1 ) ^{\anz I + 1} (\phi _1 \circledast \phi _2 )(\prod _{i \in I} a_i ) \, \bigl( \prod _{i \in I ^c } \phi _{\ge _i} (a_i )\bigr)
\end{align*}
with the convention $(\phi _1 \circledast \phi _2 )(\prod _{i \in \es} a_i ) = 1$,
the monotonic $q = 1$ case by
\begin{equation*}
(\phi _1 \lhd \phi _2 )(a_1 \ot \ldots \ot a_n ) = \phi _1 \bigl( \prod _{i , \ge _i = 1} a_i \bigr) \, \bigl( \prod _{i , \ge _i = 2} \phi _2 (a_i ) \bigr)
\end{equation*}
and the anti-monotonic $q = 1$ case by
\begin{equation*}
(\phi _1 \rhd \phi _2 )(a_1 \ot \ldots \ot a_n ) =   \bigl( \prod _{i , \ge _i = 1} \phi _1 (a_i ) \bigr) \, \phi _2 \bigl( \prod _{i , \ge _i = 2} a_i \bigr) .
\end{equation*}
Notice that the monotonic and the anti-monotonic case are not symmetric.
\par
The two parameter Boolean family is given by 
\begin{equation*}
(\phi _1 \boole _{r, s} \phi _2 )(a_1 \ot \ldots \ot a_n )  =  c(r, s, \ge )\, \phi _{\ge _1 }(a_1 ) \ldots \phi _{\ge _n } (a_n )
\end{equation*}
with
\begin{equation*}
c(r, s, \ge ) =
r^{\anz \{ i \in [n] \mit i > 1, \ge _i = 1 \} } \, s ^{\anz \{ i \in [n] \mit i > 1, \ge _i = 2 \} } .
\end{equation*}
In particular, 
\begin{align*}
(\phi _1 \boole _{r, s} \phi _2 )(a_1 a_2 ) &= s \phi _1 (a_1 ) \phi _2 (a_2 ) \\
(\phi _1 \boole _{r, s} \phi _2 )(a_2 a_1 ) &= r \phi _1 (a_1 ) \phi _2 (a_2 )  
\end{align*}
for $a_1 \in \SA _1$, $a_2 \in \SA _2$,
so that $\boole _{r, s}$ is not symmetric for $r \neq s$.
\par
If the universal product satisfies
\begin{equation*}
(\phi _1 \odot \phi _2 )(a_1 a_2 ) = (\phi _1 \odot \phi _2 )(a_2 a_1 ) = \phi _1 (a_1 ) \, \phi _2 (a_2 )
\end{equation*}
for all $a_1 \in \SA _1$, $a_2 \in \SA _2$, we say that it is normal.
It has been shown by N. Muraki \cite{Mur01,Mur02} that $\ot _1$, $\circledast _1$, $\boole _{1,1}$, $\lhd _1$ and $\rhd _1$ are the only normal u.a.u.-products in the category $\AlgProb _{1,1}$. The Muraki five are also the only positive (see Section \ref{positive})
u.a.u.-products \cite{Mur4} in $\AlgProb _{1,1}$. Moreover, $\ot _1$, $\boole _1$ and $\circledast _1$ are the only normal, symmetric u.a.u.-products and the only  positive, symmetric u.a.u.-products in $\AlgProb _{1,1}$.
}
\eexam
\bremark
{\rm
If we have a universal product in the category $\UnitalAlgProbLB$, the following procedure associates with it a universal product in the category $\AlgProbLB$.
Let
\begin{equation*}
\bigl( \SA _1 , (\SA _1  ^{(1)} , \ldots , \SA _1 ^{(m)} ) , \phi _1 \bigr), \ \bigl( \SA _2 , (\SA _2  ^{(1)} , \ldots , \SA _2 ^{(m)} ) , \phi _2 \bigr)
\end{equation*}
be objects in $\AlgProbLB$ and let $\odot _{\ei}$  be a universal product in $\UnitalAlgProbLB$. 
We form the tensor product
\begin{equation*}
\bigl( \tSA _1 \sqcup _{\ei} \tSA _2 , (\tSA _1 ^{(1)} \sqcup _{\ei} \tSA _2 ^{(1)}, \dots , \tSA _1 ^{(m)} \sqcup _{\ei} \tSA _2 ^{(m)} ) , \tphi _1 \odot _{\ei} \tphi _2 \bigr)
\end{equation*}
and put
\begin{equation*}
(\phi _1 \odot \phi _2 ) ^{(l)} :=
(\tphi _1 \odot _{\ei} \tphi _2 ) ^{(l)} \rest \SA _1 \sqcup \SA _2
\end{equation*}
which makes sense since
\begin{align*}
&\tSA _1 \sqcup _{\ei} \tSA _2 = (\tSA _1 ^{(1)} \sqcup _{\ei} \tSA _2 ^{(1)} ) \sqcup _{\ei} \ldots \sqcup _{\ei} (\tSA _1 ^{(m)} \sqcup _{\ei} \tSA _2 ^{(m)} )  \\
&   \quad \supset (\SA _1 ^{(1)} \sqcup \SA _2 ^{(1)} ) \sqcup \ldots \sqcup (\SA _1 ^{(m)} \sqcup \SA _2 ^{(m)} ) = \SA _1 \sqcup  \SA _2
\end{align*}
in the canonical way.
One checks that $\odot$ is a u.a.u.-product if this is true for $\odot _{\ei}$.
A universal product in $\AlgProbLB$ not always can be reduced to a universal product in $\UnitalAlgProbLB$. 
However, $\odot$ defines $\odot _{\ei}$ if it respects the units of the algebras.
This is, for example, the case for the tensor and the free product,  but not for the Boolean, monotone or anti-monotone product in Example \ref{five}.
}
\eremark
\bexam
\label{c-free}
{\rm
M. Bozejko and R. Speicher \cite{BSp} showed, see also \cite{BLSp}, that there is a unique u.a.u.-product $\odot _{\ei}$ in $\UnitalAlgProb _{2,1}$ such that
\begin{equation*}
(\phi _1 \odot _{\ei} \phi _2 ) ^{(2)} = \phi _1 ^{(2)}   \circledast _{\ei} \phi _2 ^{(2)}
\end{equation*}
and
\begin{align*}
& \ge \in \ma _2 , a_1 \in \SA _{\ge _1}, \ldots , a_n \in \SA _{\ge _n } , 
\phi ^{(2)} _{\ge _i } (a_i ) = 0 , i \in [n]  \\
&\quad  \implies
\bigl( (\phi _1  \odot _{\ei } \phi _2 ) ^{(1)} \bigr)(a_1 \ldots a_n )
= \phi _{\ge _1} ^{(1)} (a_1 ) \ldots \phi _{\ge _n} ^{(1)} (a_n ) .
\end{align*}
The independence coming from the u.a.u.-product $\odot$ given by the above product on $\UnitalAlgProb _{2,1}$ as in Remark \ref{c-free}, which is also called \emph{c-freeness}, gives rise to a generalized, non-commutative Brownian motion; see \cite{BSp,Ansh}.
The concept for the c-free product can be generalized to a  wider class of u.a.u.-products in $\AlgProb _{2,1}$; see \cite{Has} and Remark \ref{abstractLevy}.
}
\eexam
\bexam 
{\rm (see \cite{Voi123})
Let $\SV _1$, $\SV _2$ be vector spaces.
For $\SV _1 \sqcup \SV _2 = \bigoplus_{\ge \in \ma _2} \, \SV _{\ge}$, cf.  (\ref{freevec}), consider the four natural vector space isomorphisms
\begin{align}
\label{iso1}
\mc \oplus (\SV _1 \sqcup \SV _2 ) &\cong
(\mc \oplus \SV _1 ) \ot (\mc \oplus \bigoplus _{\mycom{\ge \in \ma _2}{\ge _1 = 2}} \, \SV _{\ge} )  \\
\nonumber
&\cong
(\mc \oplus \SV _2 ) \ot (\mc \oplus \bigoplus _{\mycom{\ge \in \ma _2}{\ge _1 = 1}} \, \SV _{\ge} ) 
\end{align}
and
\begin{align}
\label{iso2}
\mc \oplus (\SV _1 \sqcup \SV _2 ) &\cong
 (\mc \oplus \bigoplus _{\mycom{\ge \in \ma _2}{\ge _n = 2}} \, \SV _{\ge} )  
 \ot (\mc \oplus \SV _1 )\\
 \nonumber
&\cong
 (\mc \oplus \bigoplus _{\mycom{\ge \in \ma _2}{\ge _n = 1}} \, \SV _{\ge} ) 
 \ot (\mc \oplus \SV _2 ) .
\end{align}
Let $l (T)$, $T \in \End (\mc \oplus \SV _1 )$, be the element in $\End (\mc \oplus (\SV _1 \sqcup \SV _2 ))$ which equals $T \ot \id$ in the first identification of (\ref{iso1}) and let $r (T)$ be the element which equals $\id \ot T$ in the first identification of (\ref{iso2}).
Similarly, for $T \in \End (\mc \oplus \SV _2 )$ define $l(T)$ and $r (T)$ to be the elements of $\End (\mc \oplus (\SV _1 \sqcup \SV _2 ))$ using the second identifications of (\ref{iso1}) and (\ref{iso2}). 
\par
An algebraic probability space $(\SA , \phi )$, i.e. an element of $\AlgProb$, gives rise to a Gelfand-Naimark-Segal representation $\pi _{\phi} : \SA \to \End (\mc \oplus \SV _{\phi})$ which, up to isomorphisms, is uniquely determined by the properties ($\GO = (1, 0)$)
\begin{gather*}
\pi (\SA ) \GO + \mc \GO = \mc \oplus \SV _{\phi}  \\
K \subset \SV _{\phi} , \pi (\SA ) K \subset K \implies K = \{ 0\}   \\
\RP _{\GO} \pi (a) = \phi (a) \GO , \ \mbox{for all} \ a \in \SA ;
\end{gather*}
cf. \cite{Ger1,Lachs1}.
The triplet $(\pi , \SV , \GO )$ is called the \emph{GNS-triplet} of $(\SA , \phi )$.
Let  $(\SA _1 , \phi _1 )$, $(\SA _2 , \phi _2 )$ be two algebraic probability spaces with GNS-triplets 
\begin{equation*}
(\pi _1 , \SV _1 , \GO _1 ), \ (\pi _2 , \SV _2 , \GO _2 ) .
\end{equation*}
We define a representation $\pi$ of $\SA _1 \sqcup \SA _2$ on $\mc \oplus (\SV _1 \sqcup \SV _2 )$ by 
\begin{align*}
\pi (a_1 ) &= l(\pi _1 (a_1 ))   \\
\pi (a_2 ) &= l(\pi _2 (a_2 )) .
\end{align*}
Then $(\pi , \SV _1 \sqcup \SV _2, \GO )$, $\GO = (1, 0)$, is the GNS-triplet of $\phi _1 \circledast \phi _2$; see \cite{Voi85}. This GNS-triplet can equally be realized by putting
\begin{align*}
\pi (a_1 ) &= r(\pi _1 (a_1 )),   \\
\pi (a_2 ) &= r(\pi _2 (a_2 )) .
\end{align*}
It was the idea of D. Voiculescu \cite{Voi123} to involve both \lq faces\rq , i.e. the left and the right operators on $\mc \oplus (\SV _1 \sqcup \SV _2 )$, by considering objects 
\begin{equation*}
\bigl( \SA _1, (\SA _1 ^{(1)} , \SA _1 ^{(2)}) , \phi _1 \bigr), \ \bigl(\SA _2 , (\SA _2 ^{(1)} , \SA _2 ^{(2)} ), \phi _2 \bigr)
\end{equation*}
in $\AlgProb _{1,2}$ and by putting for the GNS-triplets
$(\pi _1 , \SV _1 , \GO _1 )$, $(\pi _2 , \SV _2 , \GO _2 )$ of 
$(\SA _1 , \phi _1 )$,
$(\SA _2  , \phi _2 )$
\begin{align*}
\pi (a_1 ) &= l(\pi _1 (a_1) ) \ \mbox{for } a_1 \in \SA _1 ^{(1)}  \\
\pi (a_1 ) &= r(\pi _1 (a_1) ) \ \mbox{for } a_1 \in \SA _1 ^{(2)} \\
\pi (a_2 ) &= l(\pi _2 (a_2) ) \ \mbox{for } a_2 \in \SA _2 ^{(1)}  \\
\pi (a_2 ) &= r(\pi _2 (a_2) ) \ \mbox{for } a_2 \in \SA _2 ^{(2)}
\end{align*}
to define a representation of $\SA _1 \sqcup \SA _2$ 
on $\mc \oplus (\SV _1 \sqcup \SV _2 )$.
Now put
\begin{equation*}
\phi _1 \odot \phi _2 (A) := \RP _{\GO} \bigl( \pi (A) \GO \bigr)
\end{equation*}
for $A \in \SA _1  \sqcup \SA _2 $ to obtain a u.a.u.-product $\odot$ in $\AlgProb _{1,2}$, the \emph{bi-free product} of Voiculescu.
Clearly, we have
\begin{align*}
\phi_1 \odot \phi _2 \rest           \SA _1 ^{(1)}\sqcup \SA _2 ^{(1)} &= (\phi _1 \rest \SA _1 ^{(1)} ) \circledast (\phi _2 \rest \SA _2 ^{(1)} )  \\
\phi_1 \odot \phi _2 \rest           \SA _1 ^{(2)}\sqcup \SA _2 ^{(2)} &= (\phi _1 \rest \SA _1 ^{(2)} ) \circledast (\phi _2 \rest \SA _2 ^{(2)} ) .
\end{align*}
Moreover,
\begin{align*}
\phi_1 \odot \phi _2 \rest           \SA _1 ^{(1)}\sqcup \SA _2 ^{(2)} &= (\phi _1 \rest \SA _1 ^{(1)} ) \ot (\phi _2 \rest \SA _2 ^{(2)} )  \\
\phi_1 \odot \phi _2 \rest           \SA _1 ^{(2)}\sqcup \SA _2 ^{(1)} &= (\phi _1 \rest \SA _1 ^{(2)} ) \ot (\phi _2 \rest \SA _2 ^{(1)} );
\end{align*}
see \cite{Voi123}.
Since $\odot$ is symmetric, the cumulant Lie algebras are abelian and the cumulant functional $\ln _{\odot} (\phi )$ for $ \phi \in \bigl( \RT (\SV ) \sqcup \RT (\SV )\bigr) ^* = \RT (\SV \oplus \SV ) ^*$ has the properties
\begin{align*}
\ln _{\odot} (\phi _1 \conv \phi _2 ) &= \ln _{\odot}(\phi _1 ) + \ln_{\odot} (\phi _2 )  \\
\ln _{\odot} (\phi _1 \conv \phi _2 ) \rest \i _1 \bigl( \RT (\SV )\bigr) &= \ln_{\circledast} (\phi _1 \rest \SA _1 ^{(1)} ) +  \ln_{\circledast} (\phi _2 \rest \SA _2 ^{(1)}  )\\
\ln _{\odot} (\phi _1 \conv \phi _2 ) \rest \i _2 \bigl( \RT (\SV ) \bigr) &= \ln_{\circledast} (\phi_1 \rest \SA _1 ^{(2)} ) +  \ln_{\circledast} (\phi _2 \rest \SA _2 ^{(2)}).
\end{align*}
}
\eexam
\bremark
\label{abstractLevy}
{\rm
M. Gerhold and S. Lachs \cite{Ger1,Lachs1} introduced the notion of
an \lq abstract L\'evy processes on a comonoidal system\rq . A special case is a quantum L\'evy process on a given dual semi-group with respect to a fixed independence given by a u.a.u.-product in $\AlgProb$; see \cite{Lachs1}. More generally, quantum L\'evy processes - on 
a dual semi-group and with respect to a u.a.u.-product in $\AlgProbLB$ - give rise to new examples of 
abstract L\'evy processes on a comonoidal system in the sense of \cite{Ger1,Lachs1}.
The c-free Brownian motion  of \cite{BSp,Ansh} is an example of such a L\'evy processs on a dual group of type $\RT (\SV )$ with the underlying independence given by the c-free product (Example \ref{c-free}) in $\AlgProb _{2,1}$.
See also Remark \ref{schoenberg}.
}
\eremark
\section{The central lemma}
\label{basiclemma}
The  following functors from $\AlgProbL$ to $\AlgProb$ play an important role in the classification of u.a.u.-products and for our definition of cumulants in Section \ref{cumulants}.
For $n \in \mn$ consider a functor $\SF : \AlgProbL \to \AlgProb$ such that
\begin{align}
\label{functor1}
\SF (\SA , \phi ) &= \bigl(\SA ^{\ot n} , \SF (\phi )\bigr)  \\
\label{functor2}
\SF (j) &= j^{\ot n} .
\end{align}
This means that
\begin{equation*}
\SF (\phi \circ j) (a_1 \ot \ldots \ot a_n ) = \SF (\phi ) \bigl( j (a_1 ) \ot \ldots \ot j (a_n )\bigr)
\end{equation*}
for all algebras $\SA , \SB$, algebra homomorphisms $j : \SB \to \SA$ and $\phi \in \SA ^{*d}$, $a_1 , \ldots , a_n \in \SA$.
In particular, consider the algebra $\mc \langle X_n \rangle$ and put
\begin{equation*}
F (\phi ) := \SF (\phi ) (x_1 \ot \ldots \ot x_n )
\end{equation*}
for $\phi \in \mc \langle X_n \rangle ^ {*d}$.
Then for $\phi \in \SA ^{*d}$
\begin{align*}
F (\phi \circ j (a_1 , \ldots , a _n )) &= \SF (\phi \circ j(a_1 , \ldots , a_n )) (x_1 \ot \ldots \ot x_n )  \\
 &= \SF (\phi ) (a_1 \ot \ldots \ot a_n )
 \end{align*}
 so that $\SF$ is determined by $F : \mc \langle X_n \rangle ^{*d} \to \mc $.
\par
A map
\begin{equation*}
F : \mc \lan X_n \ran ^{*d} \to \mc 
\end{equation*}
is called {\it generating} (of degree $n$ and index $d$) if the following holds. 
For each algebra $\SA$ and each element $\phi \in \SA ^{*d}$ the map
\begin{align*}
(a_1 , \ldots , a_n ) &\mapsto F\bigl( \phi \circ j (a_ 1, \ldots , a_n )\bigr)  \\
\SA ^n & \to \mc 
\end{align*}
is $n$-linear.
Complex constants are called generating of degree $0$.
It is immediate that $F$ coming from a functor $\SF$ as above is generating.
Moreover, given a generating map $F$ we put
\begin{equation*}
\SF (\phi )(a_1 \ot \ldots \ot a_n ) := F\bigl( \phi \circ j(a_1 , \ldots , a_n )\bigr)
\end{equation*}
to define a functor $\SF$ on $\AlgProbL$, which satisfies (\ref{functor1}) and (\ref{functor2}).
Indeed, the $n$-linearity makes it possible to extend $\SF (\phi )$ to a linear mapping on $\SA ^{\ot n}$.
Moreover, we have
\begin{align*}
\SF (\phi \circ j ) (a_1 \ot \ldots \ot a_n ) &= F \bigl( \phi \circ j \circ j(a_1 , \ldots , a_n )\bigr)  \\
&= F \bigl( \phi \circ j( j(a_1) , \ldots , j(a_ n )\bigr)  \\
&= \SF ( \phi ) \bigl( j(a_1 ) \ot \ldots \ot j(a_n )\bigr) .
\end{align*}
This means that we have a 1-1-correspondence between generating maps and functors $\SF$ as above.
\par
We observe that 
\begin{align}
\label{poly1}
(a_1 , \ldots , a_n ) &\mapsto \phi ^{(k_1 )} ( a_{i(1,1) } \ldots a_{i(1, s_1 )}) \ldots   \\
\nonumber
&\hspace{2cm} \ldots 
\phi ^{(k_l )} ( a_{i(l,1) } \ldots a_{i(l, s_l )})
\end{align}
is $n$-linear in $a_1 , \ldots , a_n \in \SA$, for each choice of an algebra $\SA$ and $n\in \mn$, and for all fixed $\phi \in \SA ^{*d}$ and  $(M_1 , k_1 ) \cdot \ldots \cdot (M_l , k_l )$ in the set $\STPL (n)$, $M_r = x _{i(r,1)}  \ldots x_{i(r, s_r )}$; see (\ref{labelled}) for the definition of $\STPL$.
This holds because each $a_i$, $i \in [n]$, appears exactly once as an argument in the product of the right hand side of \eqref{poly1} and because the $\phi ^{(k)}$ are linear.
Since $\phi \circ j (a_1 , \ldots , a_n )$ is as an element of $\mc \langle X_n \rangle ^{*d}$, the image of $(a_1 , \ldots , a_n )$ in (\ref{poly1})
can be written
\begin{equation*}
\S \bigl( \phi \circ j(a_1 , \ldots , a_n )\bigr) (\pi )
\end{equation*}
 with 
 \begin{equation*}
 \S \bigl( \phi \circ j(a_1 , \ldots , a_n )\bigr) \in \RS \bigl( \mc \langle  X_n \rangle ^d \bigr) ^*
 \end{equation*}
 and
\begin{equation*}
\pi = (M_1 , k_1 ) \cdot \ldots \cdot (M_l , k_l ) \in \RS \bigl( \mc \langle  X_n \rangle ^d \bigr) = \mc \bigl[ \langle X_n \rangle \times [d]\bigr].
\end{equation*}
\par
Let $P$ be a polynomial in $\mc \bigl[ \langle X_n \rangle \times [d] \bigr]$ which lies in the linear span of $\STPL (n)$, i.e.
\begin{equation*}
 P = \sum _{\pi \in \STPL (n)} \ga _{\pi}  \, \pi
 \end{equation*}
 with complex constants $\ga _{\pi}$.
Then the polynomial function
\begin{equation}
\label{poly2}
F(\phi ) = \S (\phi ) (P) = \sum_{\pi \in \STPL (n)} \ga _{\pi} \, \S (\phi ) (\pi )
\end{equation}
over $\mc (\langle X _n \rangle \times [d] ) = \mc \langle X_n \rangle ^d$ given by $P$ is generating of degree $n$ and index $d$.
Indeed, \emph{all} generating maps arise this way.
\blem
\label{mainlemma}
A map $F$ is generating of degree $n$ and index $d$ iff $F$ is a polynomial function over $\mc \langle X_n \rangle ^d$ given by a polynomial in $\mc \bigl[ \langle X_n \rangle \times [d]\bigr]$ which lies in the linear span $\mc \STPL (n)$ of $\STPL (n)$.
\elem
\pn
A proof of the remaining implication will be given at the end of this section.
\par
The right hand side of (\ref{poly2}) only depends on the values of $\phi$ on monomials with $\# (\set M ) = \vert M \vert$ that is on monomials in $\langle X_n \rangle$ \lq with no repetition of an indeterminant\rq . A first step consists in proving
\bprop
\label{polyprop1}
Let $n, d  \in \mn$ and 
let $F$ be generating of degree $n$ and index $d$.
Moreover, let $\phi , \tphi \in \mc \lan X_n \ran ^{*d} $.
\pn
Then 
\begin{align*}
&\phi (M) = \tphi (M) \mbox{ for all } M \in \lan X_n \ran \mbox{ with } \# (\set M ) = \vert M \vert  \\
&   \implies 
 F(\phi ) = F (\tphi ) .
 \end{align*}
\eprop
\pn
{\it Proof}:
We only write down the case $d = 1$, the general case can be done in exactly the same manner but with more indices.
\par
We double the number of indeterminates to obtain the algebra $\SA _1 = \mc \lan y_1, z_1, \ldots , y_n , z_n \ran$. 
Denote by $\SN$ the set of monomials 
\begin{equation*}
\xi _1 \ldots \xi _l ,  \ \xi _j \in \{ y_{i_j} , z_{i_j} \},\ i_j \in [n], \ j \in [l]
\end{equation*}
in $\SA _1$ with the following property. 
If $i_j$, $j \in [l]$, appeared before, we must have $\xi _j = z_{i_j}$, and if $i_j$ did not appear before, we must have  $\xi _j = y_{i_j}$, that is
\begin{align*}
\xi _1 &= y_{i_1}  \\
\xi _2 &= \left\{ \begin{array}{cl} y_{i_2} & \mbox{if } i_2 \neq i_1  \\
z_{i_2} & \mbox{if } i_2 = i_1 \end{array}  \right.  \\
\xi _s &= \left\{ \begin{array}{cl} y_{i_s} & \mbox{if } i_s \notin \{ i_1 , \ldots , i_{s - 1} \} \\
z_{i_s} & \mbox{if } i_s \in \{ i_1 \ldots , i_{s - 1} \} \end{array}  \right. 
\end{align*}
for $ s = 3, \ldots , l$.
For example, $y_2 y_1 y_3 z_2 y_4 z_2$ and $y_2 y_1 y_3 z_2 z_1 z_2$ are in $ \SN$, but
$y_2 y_1 y_3 z_2 y_1 z_2 \notin \SN$.
\par
Let $\psi \in \mc \langle X_n \rangle ^{*d}$.
Define the linear functional $\psi _1$ on $\SA _1$ by
\begin{equation*}
\psi _1 (\xi _1 \ldots \xi _l ) = \left\{
\begin{array}{cl}
\psi (x_{i_1} \ldots x_{i_l} )  & \mbox{if } \xi _1 \ldots \xi _l \in \SN   \\
0 & \mbox{otherwise}
\end{array}  \right.   
\end{equation*}
We claim that
\begin{align*}
\phi _1 \circ j( y_1 + z_1 , \ldots , y_n + z_n ) 
&= \phi  \\
\tphi _1 \circ j( y_1 + z_1 , \ldots , y_n + z_n ) 
&= \tphi .
\end{align*}
For $\phi$ this follows from
\begin{align*}
(\phi _1 \circ j )(x_{i_1 } \ldots x_{i_l} ) &=
\phi _1 \bigl( (y_{i_1 } + z_{i_1} ) \ldots (y_{i_l} + z_{i_l}) \bigr)   \\
&= \sum_{\mycom{\xi _1 , \ldots , \xi _l}{\xi _r \in \{ y_{i_r } , z_{i_r} \}, r \in [l]}} \phi _1 (\xi _1  \ldots \xi _l )
\end{align*}
and the fact that exactly one of the monomials $\xi _1 \ldots \xi _l$ appearing in the sum of the right hand side is in $\SN$, and that for this monomial we have $\phi _1 (\xi _1 \ldots \xi _l ) = \phi (x_{i_1} \ldots x_{i_l} )$  by the definition of $\phi _1$. Similarly for $\tphi$.
\par
Using the $n$-linearity property of the generating map $F$ twice, we have
\begin{align*}
F(\phi ) &= F \bigl(\phi _1 \circ j (y_1 + z_1 , \ldots , y_n + z_n ) \bigr)  \\
&= \sum_{\mycom{\xi _1 , \ldots , \xi _l}{\xi _r \in \{ y_{i_r } , z_{i_r} \}, r \in [l]}}
F\bigl( \phi _1 \circ j(\xi _1 , \ldots , \xi _l )\bigr)  \\
&= \sum_{\mycom{\xi _1 , \ldots , \xi _l}{\xi _r \in \{ y_{i_r } , z_{i_r} \}, r \in [l]}}
F\bigl( \tphi _1 \circ j(\xi _1 , \ldots , \xi _l )\bigr) \\
&= F(\tphi ) .
\end{align*}
Notice that equality of the sums follows because for $M \in \lan X_n \ran$
\begin{equation*}
j(\xi _1 , \ldots , \xi _n ) M \in \SN \implies
\# (\set M ) = \vert M \vert
\end{equation*}
and $\phi$ and $\tphi$ agree on monomials $M$ with $\# (\set M ) = \vert M\vert$ by assumption.\hfill$\square$
\par
We put 
\begin{equation*}
\widehat{\lan X \ran } = \{ M \in \lan X \ran \mit \# (\set M ) = \vert M \vert \}
\end{equation*}
Moreover, for $x \in X$ we put
\begin{equation*}
\lan X \ran  _x = \lan X \setminus \{ x \} \ran ; \ 
\lan X \ran ^{(x)} = \{ M \in \lan X \ran \mit  x \in \set M \} .
\end{equation*}
There is also $\widehat{\lan X \ran } _x$ and $\widehat{\lan X \ran}  ^{(x)}$ with the obvious meanings. 
We put $\widehat{\lan X_n \ran} ^ {(k)} = \widehat{\lan X_n \ran} ^ {(x_k)}$.
\bprop
\label{polyprop2}
Let $F$ be a generating map of degree $n$  and index $d$.
\pn
Then there exist
 mappings
\begin{equation*}
\ga _k ^{(M, l)} : \mc \widehat{\lan X_n \ran} _k ^{*d} \to \mc ,
\end{equation*}
$k \in [n]$, $M \in \widehat{\lan X_n \ran} ^ {(k)}$,  $l \in [d]$, 
such that 
\begin{equation*}
F(\phi ) = \sum_{(M , l) \in \widehat{\lan X_n \ran } ^{(k)} \times [d]} \ga _k ^{(M, l)} (\phi \rest \mc \widehat{\lan X_n \ran} _k ^d ) \, \phi ^{(l)}(M) .
\end{equation*}
\eprop
\pn
{\it Proof}:
Again only for $d = 1$.
We have
\begin{equation*}
\mc \wh{\lan X_n \ran} = \mc \wh{\lan X_n \ran}_k  \oplus \mc \wh{\lan X_n \ran} ^{(k)}
\end{equation*}
and
\begin{equation*}
\mc \wh{\lan X_n \ran} ^* = \mc \wh{\lan X_n \ran}_k  ^* \oplus \mc \wh{\lan X_n \ran} ^{(k)*} .
\end{equation*}
The linear functionals $\gd _M$, $M \in \wh{\lan X_n \ran} ^{(k)}$, given by $\gd _M (\phi ) = \phi (M)$, form a vector space basis of  $\mc \wh{\lan X_n \ran} ^{(k)*}$. 
The statement of the proposition means nothing else but the linearity of 
\begin{equation*}
F : \mc \wh{\lan X_n \ran} ^* = \mc \wh{\lan X_n \ran}_k  ^* \oplus \mc \wh{\lan X_n \ran} ^{(k)*}  \to \mc
\end{equation*}
in the second component.
In order to prove this linearity let $\phi , \tphi \in \hmcn ^*$ be two linear functionals which agree on $\hmcn _k$ and let $\gl \in \mc$.
Moreover, define $\theta \in \hmcn ^*$ by
\begin{equation*}
\theta (M) = \left\{ \begin{array}{cl} \phi (M) & \mbox{if } M \in \langle X_n \rangle _k  
\vspace{0.2cm}\\
\phi (M) + \gl \tphi (M) & \mbox{if } M \in \langle X_n \rangle  ^{(k)} .
\end{array} \right. 
\end{equation*}
Now consider two further indeterminates $y_k$ and $z_k$ to double the indeterminate $x_k$ and put
\begin{equation*}
(\xi _1 , \dots , \xi _{n + 1} ) = (x_1 , \ldots , x_{k - 1}, y_k , z_k , x_{k + 1} , \ldots , x_n )
\end{equation*}
and
\begin{equation*}
\SA _1 := \mc \langle \xi _1, \ldots \xi _{n + 1} \rangle .
\end{equation*}
Next,  define $\ttheta \in (\SA _1 ) ^*$ by ($\xi = \xi _{i_1} \ldots \xi _{i_l}$ a monomial in $\SA _1$)
\begin{equation*}
\ttheta (\xi ) = \left\{
\begin{array}{cl}
0 & \mbox{if } k, k + 1 \in \{ i_1 , \ldots , i_l \} 
\\
\phi (x_{i_1} \ldots x_{i_l} ) & \mbox{if } k + 1 \notin \{ i_1 , \ldots , i_l \} 
\\
\tphi (x_{i_1} \ldots x_{i_l} ) & \mbox{if } k  \notin \{ i_1 , \ldots , i_l \} , \
k + 1  \in \{ i_1 , \ldots , i_l \} 
\end{array}
\right.
\end{equation*}
and consider the algebra homomorphism $j : \mcn \to \SA _1$ which is given by
\begin{equation*}
j (x_ i ) = 
\left\{ 
\begin{array}{cl}
\xi _i & \mbox{if } i \neq k  \\
x_k + \gl z_k & \mbox{if } i = k  .
\end{array}
\right. 
\end{equation*}
We have $\ttheta \circ j = \theta$.
Indeed, if $M = x_{i_1} \ldots x_{i_l}$ is in $\langle X_n \rangle$ with $k \notin \{i_1 , \ldots , i_l \}$, then 
\begin{equation*}
(\ttheta \circ j )(M) = \phi (M) = \tphi (M) = \theta (M)
\end{equation*}
whereas in the other case $k \in \{i_1 , \ldots , i_l \}$
\begin{align*}
(\ttheta \circ j)(M)   
&= \ttheta (x_{i_1} \ldots x_{i_{k - 1}} (y_k + \gl z_k ) x_{i_{k + 1}} \ldots x_{i_l} )  
\\
&= \ttheta (x_{i_1} \ldots x_{i_{k - 1}} y_k  x_{i_{k + 1}} \ldots x_{i_l} ) \\
& \hspace{1cm} +
 \gl \ttheta (x_{i_1} \ldots x_{i_{k - 1}} z_k  x_{i_{k + 1}} \ldots x_{i_l} )
 \\
&= \phi (M) + \gl \tphi (M)  \\
&= \theta (M) .
 \end{align*}
Using that $F$ is generating, we end up with
\begin{align*}
F(\theta ) &=  F (\ttheta \circ j ) \\
&= F \bigl(\ttheta \circ j(x_1 , \ldots , x_{k - 1}, y_k + \gl z_k , x_{k + 1} , \ldots , x_n ) \bigr) \\
&= F \bigl(\ttheta \circ j(x_1 , \ldots , x_{k - 1}, y_k , x_{k + 1} , \ldots , x_n ) \bigr) \\
 &\hspace{1cm} +  \gl F \bigl(\ttheta \circ j(x_1 , \ldots , x_{k - 1}, y_k , x_{k + 1} , \ldots , x_n ) \bigr)  \\
 &= F (\phi ) + \gl F (\tphi )
 \end{align*}
 which proves linearity of $F$ in the second component.\hfill$\square$
\par
For simplicity, from now on we will drop the hats in $\wh{\lan X_n \ran} $ and $\wh{\lan X_n \ran }_k$ etcetera
 and bear in mind that the values of our functionals only depend on monomials without repetitions of indeterminates.
\bprop
\label{propnoch}
{\rm \bf (a)} 
Let $n \in \mn$ and $(M, l)  \in {\lan X_n \ran} \times [d]$, $\vert M \vert \leq n$,  and let $k \in [n]$ with $x_k \in \set M$.
\pn
If $\phi , \psi \in \mc {\lan X_n \ran} _k ^{*d} $ agree on $\mc {\lan X_n \setminus \set M \ran} ^d$, then we have
\begin{equation*}
\ga ^{(M, l)} _k (\phi ) = \ga ^{(M, l)} _k (\psi ) .
\end{equation*}
The functions $\ga ^{(M, l)} _k : \mc {\lan X_n \setminus \set M \ran} ^{*d }\to \mc $ are the same for all $k \in \set M$.
\pn
{\rm  \bf (b)} 
In particular, if
$\vert M \vert = n $, we have that
$\ga ^{(M, l)} _k$ is constant and does not depend on $k$.
\eprop
\pn
{\it Proof}: Only for $d = 1$.
\pn
{\rm \bf (a)} :
Let $k, p \in [n]$ with $x_k , x_p \in \set M$. By Proposition \ref{polyprop2} we have for $\phi \in \mcn ^*$
\begin{align*}
F (\phi ) &= \sum_{K \in \lan X_n \ran ^{(k)}} \ga _k ^{(K)} (\phi \rest \mc \lan X_n \ran _k ) \phi (K)   \\
&= \sum_{P \in \lan X_n \ran ^{(p)}} \ga _p ^{(P)} (\phi \rest \mc \lan X_n \ran _p ) \phi (P) .
\end{align*}
Splitting each of the two sums we have
\begin{align*}
&\ga _k ^{(M)} (\phi \rest \mc \lan X_n \ran _k )\phi (M) +
\sum_{\mycom{K \in \lan X_n \ran ^{(k)}}{K \neq M}} \ga _k ^{(K)} (\phi \rest \mc \lan X_n \ran _k ) \phi (K)   \\
&= 
\ga _p ^{(M)} (\phi \rest \mc \lan X_n \ran _p ) \phi (M) + \sum_{\mycom{P \in \lan X_n \ran ^{(p)}}{P \neq M}} \ga _p ^{(P)} (\phi \rest \mc \lan X_n \ran _p ) \phi (P) 
\end{align*}
and
\begin{align}
\label{trick}
& \bigl( \ga _k ^{(M)} (\phi \rest \mc \lan X_n \ran _k ) - \ga _p ^{(M)} (\phi \rest \mc \lan X_n \ran _p )  \bigr) \phi (M)
   \\
   \nonumber
&= 
  \sum_{\mycom{P \in \lan X_n \ran ^{(p)}}{P \neq M}} \ga _p ^{(P)} (\phi \rest \mc \lan X_n \ran _p ) \phi (P)  \\
  \nonumber
 & \hspace{1,5cm} -  \sum_{\mycom{K \in \lan X_n \ran ^{(k)}}{K \neq M}} \ga _k ^{(K)} (\phi \rest \mc \lan X_n \ran _k ) \phi (K) .
\end{align}
Since $x_k, x_l \in \set M$ we have $M \notin \lan X_n \ran _l \cup \lan X_n \ran _k$, which means that the right hand side of (\ref{trick}) does not depend on $\phi (M)$ for each choice of $\phi \in \mc \lan X_n \ran ^*$.
From this it follows that
\begin{equation}
\label{kp}
\ga _k ^{(M)} (\phi \rest \mc \lan X_n \ran _k ) = \ga _p ^{(M)} (\phi \rest \mc \lan X_n \ran _p )
\end{equation}
for all $\phi \in \mc \lan X_n \ran ^*$.
\par
Now let  $\set M = \{ x_{i_1} \ldots x_{i_l} \}$, $\vert M \vert = l$, $1 \leq l < n$  and $\phi \in \mcn ^*$. 
We will show that
\begin{equation*}
\ga _{i_1} ^{(M)} (\phi  \rest \mc \lan X_n \ran _{i_1} ) =
\ga _{i_1} ^{(M)} (\psi \rest \mc \lan X_n \ran _{i_1} )
\end{equation*}
where $\psi$ is the linear functional on $\mcn$ which agrees with $\phi$ on elements of
$\mc \lan X_n \setminus \set M \ran$ and is 0 on other elements of $\mcn$.
Change $\phi$ to $\phi _1 \in \mcn ^*$ by putting
\begin{equation*}
\phi _1 (W) = 
\left\{ \begin{array}{cl} \phi (W) & \mbox{if } W \in \lan X_n \ran _{i_1}  \\
0 & \mbox{otherwise} .
\end{array}
\right. 
\end{equation*}
Then, of course, 
\begin{equation*}
\ga _{i_1} ^{(M)} (\phi \rest \mc \lan X_n \ran _{i_1} ) = \ga _{i_1} ^{(M)} (\phi _1 \rest \mc \lan X_n \ran _{i_1} ) .
\end{equation*}
Next define $\phi _2 \in \mcn ^*$ by
\begin{equation*}
\phi _2 (W) = 
\left\{ \begin{array}{cl} \phi _1 (W) & \mbox{if } W \in \lan X_n \ran _{i_1} \cap \lan X_n \ran _{i_2}  \\
0 & \mbox{otherwise} .
\end{array}
\right.
\end{equation*}
By (\ref{kp}) and the fact that $\phi _2 \rest \mcn _{i_2} = \phi _1 \rest \mcn _{i_2}$ we have
\begin{align*}
\ga _{i_1} ^{(M)} (\phi _2 \rest \mc \lan X_n \ran _{i_1} )
&= \ga _{i_2} ^{(M)} (\phi _2 \rest \mc \lan X_n \ran _{i_2} )  \\
&= \ga _{i_2} ^{(M)} (\phi _1 \rest \mc \lan X_n \ran _{i_2} )  \\
&= \ga _{i_1} ^{(M)} (\phi _1 \rest \mc \lan X_n \ran _{i_1} )  .
\end{align*}
Go on by defining $\phi _3,   \ldots \phi _{l - 1}$ and finally $\phi _l = \psi$ by
\begin{equation*}
\phi _l (W) = 
\left\{ \begin{array}{cl} \phi _{l - 1} (W) & \mbox{if } W \in \lan X_n \ran _{i_1} \cap \ldots \cap \lan X_n \ran _{i_l}  \\
0 & \mbox{otherwise}
\end{array}
\right.
\end{equation*}
to obtain
\begin{align*}
\ga _{i_1} ^{(M)} (\phi _l \rest \mc \lan X_n \ran _{i_1} )
&= \ga _{i_1} ^{(M)} (\phi _{l - 1} \rest \mc \lan X_n \ran _{i_1} )  \\
& \hspace{0.2cm} \vdots  \\
&= \ga _{i_1} ^{(M)} (\phi _1 \rest \mc \lan X_n \ran _{i_1} )  \\
&= \ga _{i_1} ^{(M)} (\phi \rest \mc \lan X_n \ran _{i_1} ) .
\end{align*}
This is the first statement of {\bf (a)} if we choose $i_1 = k$, the second following directly from (\ref{kp}).\hfill$\square$
\bn
Put $\ga ^{(M, l)} := \ga _k ^{(M,l)}$, $k \in \set M$.
\bprop
\label{induct}
Let $(M, l) \in \langle X_n \rangle  \times [d], n \in \mn$. 
\pn
Then the mapping 
\begin{equation*}
\ga ^{(M, l)} : \mc \langle X_n \setminus \set M \rangle ^{*d} \to \mc 
\end{equation*}
(where $\ga ^{(M; l)} \in \mc$ if $\vert M \vert = n$)
is generating of degree $n - \vert M\vert$ and index $d$.
\eprop
\pn
{\it Proof}: ($d = 1$)
\pn
Let $X _n \setminus \set M = \{ x_{i_1} , \ldots , x_{i_r} \}$ and 
let $(\SA , \phi ) $ be an object of $\AlgProb$.
Moreover, let $a_{i_1} , \ldots , a_{i_r}  \in \SA$.
We must show that
\begin{equation*}
(a_{i_1} , \ldots , a_{i_r} ) \mapsto \ga ^{(M)} \bigl(\phi \circ j(a_{i_1}, \ldots a_{i_r})\bigr)
\end{equation*}
is $r$-linear.
\par
We have
\begin{equation*}
\SA \sqcup \mcn = \SA \, \oplus \, \mcn \, \oplus \,\bigoplus_{\mycom{\ge \in \ma _2}{\vert \ge \vert \geq 2}} \SA _{\ge} 
\end{equation*}
 with $\SA _1 = \SA , \ \SA _2 = \mcn$.
Let $\tphi$ be the linear functional on $\SA \sqcup \mcn$ which equals $\phi$ on $\SA$, $\gd _M$ on $\mcn$ and which is $0$ on $\bigoplus_{\vert \ge \vert \geq 2} \SA _{\ge}$.
Denote by $j : \mcn \to \SA \sqcup \mcn$ the algebra homomorphism with
\begin{equation*}
j(x_ i ) = \left\{ 
\begin{array}{cl}
a_i & \mbox{if } i \in X_n \setminus \set M   \\
 0 &  \mbox{otherwise} .
\end{array}
\right. 
\end{equation*}
Let $k \in \set M$.
By Propositions \ref{polyprop2} and \ref{propnoch}
\begin{align}
\label{poly4}
&F(\tphi \circ j ) \\
\nonumber
&= \sum_{N \in \lan X_n \ran ^{(k)}} \ga ^{(N)} \bigl( \tphi \circ j \rest \mc \lan X_n \setminus \set N \ran \bigr) (\tphi \circ j )(N) .
\end{align}
Moreover, we have that 
\begin{equation*}
\tphi \circ j \rest \mc \lan X_n \setminus \set M \ran = \phi \circ j(a_{i_1} , \ldots , a_{i_r})
\end{equation*}
and $(\tphi \circ j )(N) = \gd _M (N)$ for $N \in \lan X_n \ran ^{(k)}$.
Equation (\ref{poly4}) becomes
\begin{equation*}
F(\tphi \circ j ) = \ga ^{(M)} \bigl(\phi \circ  j(a_{i_1} , \ldots , a_{i_r})\bigr) .
\end{equation*}
The left side of this equation is $r$-linear in $a_{i_1}, \ldots a_{i_r}$ because $F$ is generating of degree $n$.\hfill$\square$
\bn
{\it Proof} of Lemma \ref{mainlemma}:
\pn
For $n = 1$ we have by Proposition \ref{polyprop2}
\begin{equation*}
F (\phi ) = \ga ^{(x, 1)} \phi ^{(1)} (x) + \ldots +  \ga ^{(x, d)} \phi ^{(d)} (x)
\end{equation*}
which is (\ref{poly2}) with $\ga _{\{(x , l)\}}  = \ga ^{(x, l)}$.
To prove the lemma for $n + 1$ under the assumption that it holds for natural numbers $\leq n$, we use again Proposition \ref{polyprop2} to write (with the obvious interpretation if $\vert M\vert = n + 1$)
\begin{equation}
F (\phi ) = \sum_{ (M, l ) \in \langle X_{n + 1} \rangle ^{(1)} \times [d]}
\ga  ^{(M, l )} \bigl( \phi \rest \mc \langle X_{n + 1} \setminus \set M \rangle ^d \bigr) \,
\phi ^{(l)} (M)
\end{equation}
which by the induction hypothesis (\ref{induct}) equals
\begin{align}
\label{hypo}
\sum _{(M, l)} \Bigl  (\sum _{\pi} \ga _{\pi} ^{(M, l)}  
\prod_{ (N, r ) \in \pi} \phi ^{(r)} (N)\Bigr) \phi ^{(l)} (M)
\end{align}
where the first sum is over all $(M, l) $ in $\langle X_{n + 1} \rangle ^{(1)} \times [d]$ and the second over all $\pi$ in $\STPL ( X_{n + 1} \setminus \set M)$, and where 
$\ga _{\pi} ^{(M, l)}$ are suitable complex numbers.
By $(N, r ) \in \pi$ we mean that $(N, r)$ appears in $\pi$ as a factor.
The product in (\ref{hypo}) is nothing else but $\S (\phi )(\pi )$.
Each  $\pi$ in   $\STPL  (n + 1)$ has a unique factor 
$(M_{\pi }, l _{\pi } )$ such that $x_1 \in \set M_{\pi}$.
Therefore, the above expression (\ref{hypo}) is equal to
\begin{align*}
\sum_{\pi \in \STPL  (n + 1)} \ga _{\pi \setminus \{ (M_{\pi}), l _{\pi})\}} ^{(M_{\pi} , l _{\pi})} \prod_{ (M, l) \in \pi } \phi ^{(l)} (M)
\end{align*}
which gives (\ref{poly2})  for degree $n + 1$.
\hfill $\square$
\par
Finally, we reformulate Lemma \ref{mainlemma}.
Let $\SV$ be a vector space and $\SS (\SV )$ the symmetric tensor algebra over $\SV$; see Section \ref{pres}. 
\bprop
\label{mainlemma2}
For a generating map $F$  of degree $n$ and index $d$ there exist unique mappings
\begin{equation*}
\gs _{\SA} : \SA ^n \to \RS (\SA ^d ) = \RS (\SA ) ^{\ot d},  \ \SA \mbox{ an algebra}
\end{equation*}
such that
\begin{equation*}
F \bigl( \phi \circ j(a_1 , \ldots , a_n )\bigr)  = \bigl( \RS (\phi ) \circ \gs _{\SA }  \bigr) (a_1 , \ldots , a_n )
\end{equation*}
for all $a_1 , \ldots , a_n \in \SA$, $\phi \in \SA ^{*d}$.
\eprop
\pn
{\it Proof}:
Uniqueness is clear, because a polynomial over the field of complex numbers is uniquely determined by its polynomial function and we have that $\gs _{\SA } (a_1, \ldots , a_n )$ is a complex polynomial of the polynomial algebra $\RS (\SA ^d)$.
\par
Given $F$ put ($M(a_1 , \ldots , a_n ) := j (a_1 , \ldots , a_n )(M)$)
\begin{equation*}
\gs _{\SA } (a_1 , \ldots , a_n ) = 
\sum _{\pi \in \STPL (n)}
\ga _{\pi} \, 
\prod_{ (M, l ) \in \pi} M(a_1 , \ldots a_n )^{(l)}
\end{equation*}
where the product $\prod$ is taken in $\RS (\SA ^d )$ and for $l \in [d]$, $a \in \SA$, we denote by $a ^{(l)}$ the embedding of $a$ into the $l$-th coordinate in $\SA ^d \subset \RS ( \SA ^d )$.\hfill$\square$
\par
The algebra $\mcn$ is an $\mn $-graded algebra if we put the grade of a polynomial $m \in \mcn$ equal to $\vert M\vert$.
The symmetric algebra $\S (\SV )$ over an $\mn$-graded vector space $\SV$ is an $\mn _0$-graded algebra in the natural manner.
If we equip $\mcn ^d$ with the $\mn$-grading of the direct sum, the algebra $\S (\mcn ^d )$ becomes an $\mn _0$-graded algebra. We reformulate \ref{mainlemma} to obtain
\bthm
There is a bijection $P \mapsto F_P$ between  polynomials $P$ in $\mc \STPL (n)$  and generating maps $F$  of degree $n$ and index $d$, which is given by
\begin{align*}
F_P (\phi ) &= \S (\phi ) (P ), \ \phi \in \mcn ^{*d} \\
P &= \gs _{\mc \langle X_n \rangle}(x_1 , \ldots , x_n ) .
\end{align*}
\ethm
\pn

\section{Generating families and positivity}
\label{positive}
A \emph{generating family} of index $d$ over a set $I$ is a family 
$\bF = (F_ {\ge})_{\ge \in \ma (I)}$ where $F_{\ge }$ is a generating map of degree $n = \vert \ge \vert$ and index $d$. In the case $I = [k]$, $k \in \mn$, we will speak of a $k$-fold generating family. For an algebra $\SA$ and $\phi \in \SA ^{*d}$ we let $\bF (\phi )$ be the linear functional on $ \SA ^{\sqcup I}$ (see (\ref{freeprod})) given by ($a_1 \ldots a_n \in \SA _{\ge}$)
\begin{align*}
\bF (\phi ) (a_1 \ldots a_m ) = F_{\ge } \bigl( \phi \circ j(a_1 , \ldots , a_m )\bigr) .
\end{align*}
\bremark
{\rm
A generating family is a functor $\SF : \AlgProb _d \to \AlgProb$ such that
\begin{align*}
\SF (\SA , \phi ) &= \bigl(\SA ^{\sqcup I} , \SF (\phi )\bigr)  \\
\SF (j) &= j^{\sqcup I} ;
\end{align*}
cf. (\ref{functor1}) and (\ref{functor2}).
}
\eremark
\pn
We say that the family $\bF$ is {\it positive} if the following holds. For each $*$-algebra $\SA$ and strongly positive $\phi$ (see Section \ref{pres}) on $\SA ^d$ we have that $\bF (\phi )$ is positive on the $*$-algebra $\SA ^{\sqcup I}$.
\par 
Denote by $\STOPL (n)$ the subset of $\STPL (n) $ of  \lq right\rq \ ordered elements of $\STPL (n)$, i.e. of elements 
\begin{equation*}
(M_1 , k_1 )\cdot \ldots \cdot ( M_l , k_l ), \ M_r = x _{i(r,1)}  \ldots x_{i(r, s_r )}
\end{equation*}
with $i(r, 1 ) < \ldots < i(r, s_r )$, $r \in [l]$.
Elements of $\STPL (n) \setminus \STOPL (n)$ are called \emph{wrong-ordered}.
\bremark
{\rm The result of the following theorem in the positive case is due to N. Muraki; see \cite{Mur4}.
Notice that the vanishing of the \lq wrong-ordered coefficients\rq \ $\ga _{\pi} ^{(\ge )}$, $\pi \notin \STOPL$, follows from positivity without any further assumptions, e.g. associativity conditions.}
\eremark
\bthm
\label{polfam}
Let
$\bF$ be a generating family of index $d$ over a set $I$ and let $\ge \in \ma (I)$.
\pn
{\bf (a)} 
Then
there are uniquely determined constants $\ga _{\pi} ^{(\ge )} \in \mc $,
$\pi \in \STPL (n)$, such that
\begin{align}
\label{pos1}
& \bF (\phi ) (a_1 \ldots a_n ) \\
& \ = \sum_{\pi \in \STPL (n)} \ga _{\pi} ^{(\ge )} \, 
\prod_{ (M, l ) \in \pi} \phi ^{(l)} \bigl( j (a_1, \ldots , a_n ) (M)\bigr)
\nonumber
\end{align}
for all $\phi \in \SA ^{*d}$, $a_1 \ldots a_n \in \SA _{\ge }$. 
\pn
{\bf (b)} If $\bF$ is positive, then there are uniquely determined constants 
\pn
$\ga _{\pi} ^{(\ge )} \in \mc $, $\pi \in \STOPL (n)$, such that
\begin{align*}
& \bF (\phi ) (a_1 \ldots a_n ) \\
& \ = \sum_{\pi \in \STOPL (m)} \ga _{\pi} ^{(\ge )} \, 
\prod_{ (M, l ) \in \pi} \phi ^{(l)} \bigl( j (a_1, \ldots , a_n ) (M)\bigr)
\end{align*}
for all $\phi \in \SA ^{*d}$, $a_1 \ldots a_n \in \SA _{\ge }$.  
\ethm
\pn
{\it Proof}: { (a)}  follows from Lemma \ref{mainlemma}, since $F_{\ge}$ is 
generating of degree $n = \vert \ge \vert$. 
\par
For { (b)}  we give a sketch of the proof in the case $d = 1$, which  copies the argument in \cite{Mur4}.
Let $\bF$ be a positive generating family.
Then we have the formula for $\bF$ given in {\bf (a)}  with constants $\ga _{\pi} ^{(\ge )}$, $\ge \in \ma _I$, $\pi \in \STP (n)$. 
We must show that positivity implies $\ga _{\pi} ^{(\ge )} = 0$ for $\pi \notin \STOP (n)$.
\par
Let $\tau \in \ma (I)$ , $\vert \tau \vert = n$ and let $\eta = M_1 ,\cdot\ldots \cdot M_r \notin \STOP (n)$. Then we may assume that $M_1$ is wrong ordered.
Consider the Hilbert space
\begin{equation*}
H = \bigoplus_{s = 1} ^r \mc ^{d _s}; \ d_s = \vert M_s \vert
\end{equation*}
and denote by $\{ e ^{(s)} _1, \ldots , e^{(s)} _{d_s} \}$ the canonical orthonormal basis of $\mc ^{d_s}$ as a subspace of $H$. 
Put $\GO = e ^{(1)} _1 + \ldots + e^{(r)} _1$.
Define a representation $\rho$ on $H$ of the $*$-algebra 
\begin{equation*}
\mc \lan n, * \ran :=  \mc \lan x_1 , \ldots , x_n , x _1 ^* , \ldots , x_n ^* \rangle
\end{equation*}
generated by $x_1 , \ldots , x_n$  by
putting $\rho (x_i ) e ^{(s)} _l$ equal to $0$ unless $x_i \in \set M_s$ \emph{and} $x_i$ equals the 
$(l - 1)$-st factor in the monomial $M_s$.
In the latter case we put $\rho (x_i ) e ^{(s)} _l =  e ^{(s)} _{l - 1}$.
Here $l$ and $l - 1$ are to be understood modulo $d_s$.
Let $\phi$ be the strongly positive linear functional on $\mc \lan n , * \ran$ given by
\begin{equation*}
\phi (T) = \lan \GO , \rho (T ) \GO \rangle , \ T \in \mc  \lan n , * \ran .
\end{equation*}  
One checks that
\begin{equation}
\label{pos2}
 \sum_{\pi \in \STP (n)} \ga _{\pi } ^{(\tau )} \, 
\prod_{ M \in \pi} \phi  ( M) = \ga ^{(\tau )} _{\eta} .
\end{equation}
It follows from ({\ref{pos1}) with $\SA = \mc \lan n , * \ran$, $\ge = \tau$, $a_i = x_i$, $i \in [n]$, and from (\ref{pos2}) that we have
\begin{equation*}
\ga ^{(\tau )}_{\eta } = F _{\tau} (\phi )(x_1 \ldots x_n) .
\end{equation*}
Since $M_1$ is wrong ordered there are $p, q \in [n]$, $p < q$, such that 
\begin{equation*}
M_1 = x_{i(1)} \ldots x_q \ldots x_p \ldots x_{i({d_1})}
\end{equation*}
Now we use that $\bF$ is positive. Since $\phi$ strongly positive 
on $\mc \lan n , * \ran$, we have that $\bF (\phi )$ is positive on $(\mc \lan n , * \ran )^{\sqcup I}$.
By Cauchy-Schwartz 
\begin{align}
\label{pos4}
&\vert \bF (\phi )(x_1 \ldots x_p \ldots x_q \ldots x_n ) \vert ^2  \\
&\hspace{0.5cm} \leq \bF (\phi )\bigl( (x_1 \ldots x_p )(x_1 \ldots x_p )^* \bigr) \nonumber \\
&\hspace{1.5cm} \bF (\phi ) \bigl( (x_{p + 1} \ldots x_q \ldots x_n ) ^* (x_{p + 1} \ldots x_q \ldots x_n ) \bigr)
\nonumber
\end{align}
Next, again by (\ref{pos1}), this time  for 
\begin{gather*}
\SA = \mc \lan n , * \ran , \ \ge = (\tau _1 , \dots , \tau _p , \tau _{p - 1} , \ldots , \tau _1 ), \ n = 2p -1  , \\
(a_1 , \ldots, a_{2p - 1} ) = (x_1 , \ldots , x_{p - 1} , x_p  x_p  ^*, x_{p - 1}^* , \ldots , x_1 ^* ),
\end{gather*}
we have
\begin{align}
\label{pos3}
&F (\phi ) (x_1 \ldots x_p \, x_p ^* \ldots x_1 ^* )  \\
&\hspace{0.3cm} = 
\sum_{\pi \in \STP (2p - 1)} \ga _{\pi} ^{(\ge )} \prod_{M \in \pi} 
\phi \bigl( j(x_1, \ldots , x_p x_p ^* , \ldots , x_1 ^* )(M) \bigr) .
\nonumber
\end{align}
Let $\pi \in \STP (2p - 1)$. Then there is a unique $M \in \pi$ with $p \in \set M$.
Since $q > p$, the indeterminate $x_q$ cannot appear as a factor of the monomial 
\begin{equation*}
K  = j(x_1, \ldots , x_p x_p ^* , \ldots , x_1 ^* )(M).
\end{equation*}
One checks that 
\begin{equation*}
\lan \GO , \rho (K) \GO \ran = \lan e ^{(1)} _1 , \rho (K) e ^{(1)} _1 \ran = 0
\end{equation*}
which means that (\ref{pos3}) and thus (\ref{pos4}) is 0, from which it follows that the coefficient of the wrong ordered partition $\eta$ is $0$. \hfill$\square$
\par
We apply Theorem \ref{polfam} to universal products.
For $\ge \in \ma (I)$, $i \in I$, put ($n = \vert \ge \vert$)
\begin{equation*}
X ^{(i)} = \{ x_l  \in X_n \mit\ge _l = i \} .
\end{equation*}
We associate with  a universal product in $\AlgProbL$ and $k \in [d]$ a generating family $\bF ^{(k)} = (F _{\ge} ^{(k)} )_{\ge \in \ma (I)}$ of  index $d$ over $I$ by putting
\begin{equation*}
F_{\ge} ^{(k)}(\phi )  : = \bigl( \odot_{i \in I} \, \phi _i \bigr) ^{(k)} (\iota _{\ge _1 } (x_1 ) \ldots \iota _{\ge _n } (x_n ))
\end{equation*}
where $\iota _i : \mc \lan X_n \ran \to \mc \lan X_n \ran ^{\sqcup I}$ is the natural embedding of an $i$-th copy of $\mc \lan X_n \ran$ and where $\phi _i : \mc \lan X ^{(i)} \ran ^* \to \mc$ is the restriction of $\phi \in \mc \lan X_n \ran ^*$.
\par
It follows from Theorem \ref{polfam} and the fact that 
\begin{equation*}
F_{\ge} (\phi ) = \bigl( F_{\ge} ^{(1)} (\phi ), \ldots, F_{\ge} ^{(d)}(\phi )\bigr)
\end{equation*}
does not depend on the values of $\phi$ on monomials from mixtures of $\mc \lan X ^{(i)} \ran$, $i \in I$, that we have the following result which we only write down for the case $I = [k]$.
\bthm
\label{coeff}
Let 
$\odot$ be a universal product in $\AlgProbL$ over $I = [k]$, $k \in \mn$, and
let $\ge \in \ma _k$.
Moreover, let $(\SA _1 , \phi _1 ), \ldots , (\SA _k , \phi _k )$  be  $k$ objects in $\AlgProbL$,  and let $a_1 \ldots a_n \in \SA _{\ge}$.
\pn
Then  there exist uniquely (by $\odot$) determined constants 
$\ga ^{(\ge )} _{\pi _1 , \ldots \pi _k } \in \mc ^d$, $\pi _i \in \STPL (X^{(i)})$, such that
\begin{align}
\label{formulaproducts}
&\bigl( \odot _{i = 1} ^k \, \phi _i \bigr) (a_1 \ldots a_n )  \\
\nonumber
& \  = \sum_{\pi _1 \in \STPL (X ^{(1)} )} \ldots   \sum_{\pi _k \in \STPL (X^{(k)} )} 
\ga ^{(\ge )} _{\pi _1 , \ldots \pi _k } \\
\nonumber
&\hspace{1cm}  \prod_{ (M _1 , l _1 ) \in \pi _1} \phi _1 ^{(l_1 )} \bigl( j (a_1, \ldots , a_n ) (M_1 )\bigr)
\ldots   \\
\nonumber
&\hspace{3cm} \ldots
\prod_{ (M _k , l _ k ) \in \pi _k} \phi _k ^{(l_k )} \bigl( j (a_1, \ldots , a_n ) (M_k )\bigr) .
\end{align}
 \hfill $\square$
\ethm
\bremark
{\rm
The version of Theorem \ref{coeff} for a universal product $\odot$ in $\AlgProbLB$ formally looks the same,
but $\ma _k$ has to be replaced by $\ma ([k] \times [m])$ and for $\ge = (\ge _1 , \ldots , \ge _n )$, $\ge _i = (\ge _{i,1}, \ge _{i,2} ) \in [k] \times [m]$, $i \in [k]$, we put 
\begin{equation*}
X ^{(i)} = \{ x_l \in X_n \mit \ge _{l,1} = i \} .
\end{equation*}
}
\eremark
\bremark
{\rm
A universal product is called \emph{positive} if the corresponding generating families are positive. It follows from Theorem \ref{polfam} (b) that, for positive universal products, in Theorem \ref{coeff} $\STPL$ can be replaced by $\STOPL$.
}
\eremark
\bremark
\label{schoenberg}
{\rm
It is an open problem to prove the \lq Schoenberg correspondence\rq , cf. \cite{SchVo} for the case $d = m = 1$, for positive u.a.u-products in the general $\AlgProbLB$ case, i.e. the correspondence between conditionally positive $d$-tuples of linear functionals on $D ^{\sqcup m}$, $D$ a dual semi-group with involution, and convolution semi-groups of $d$-tuples of states on $D^{\sqcup m}$.
It should suffice to show the correspondence in the cases $D = \RT (\SV )$ by constructing a representation of the quantum L\'evy process on a suitable \lq Fock space\rq, the general case then following from an approximation argument like in \cite{SchVo}.
}
\eremark

\section{The Lachs functor}

Using the central Lemma \ref{mainlemma}, we will \lq reduce\rq \ a pair $(\SD , \odot )$, $\SD$ a dual group, $\odot$ a u.a.u.-product, to a commutative bialgebra.
\bprop
Let
 $\odot$ be a universal product in $\AlgProbLB$.
 \pn
 Then there exist unique mappings
\begin{equation*}
\gs _{\SA _1 , \SA _2}  : (\SA _1 \sqcup \SA _2 ) ^d \to
\RS (\SA _1 ) ^{\ot d} \ot \RS (\SA _2 ) ^{\ot d}
\end{equation*}
such that
\begin{equation}
\label{formula3}
(\phi _1 \odot \phi _2 ) = \bigl( \RS (\phi _1 ) \ot \RS (\phi _2 ) \bigr) \circ \gs _{\SA _1 , \SA _2} 
\end{equation}
\eprop
\pn
{\it Proof}: Let $\ge \in \ma ([2] \times [m] )$, $a_1 \ldots a_n \in \SA _{\ge}$.
We put
\begin{align}
\label{formula2}
&\gs _{\SA _1 , \SA _2} (a_1 \ldots a_n ) \\
\nonumber
&= \sum_{\pi _1 \in \STPL (X ^{(1)} )} \, \sum_{\pi _2 \in \STPL (X^{(2)} )} 
\ga ^{(\ge )} _{\pi _1 ,  \pi _2 } \\
\nonumber
&\hspace{1cm} \bigl( \prod_{ (M _1 , l _1 ) \in \pi _1} \, M_1 (a_1, \ldots , a_n ) ^{(l_1)} \bigr) \ot
\bigl( \prod_{ (M _2 , l _ 2 ) \in \pi _2} \, M_2 (a_1 , \ldots , a_n )^{(l_2 )}  \bigr) .
\end{align}
with the constants $\ga ^{(\ge )} _{\pi _1 ,  \pi _2 } \in \mc ^d$ of (\ref{formulaproducts}) of Theorem \ref{coeff} ($k = 2$).
The right hand side of (\ref{formula2}) is a $d$-tuple of elements of $\RS (\SA _1 ) ^{\ot d} \ot \RS (\SA _2 ) ^{\ot d }$ so that (\ref{formula2}) defines a mapping from $(\SA _1 \sqcup \SA _2 )^d$ to $\RS (\SA _1 ) ^{\ot d} \ot \RS (\SA _2 ) ^{\ot d }$. Equation (\ref{formula3}) is (\ref{formulaproducts}).
\hfill $\square$
\bprop
\label{axioms-sigma}
{\bf (a)}
The mappings $\gs _{\SA _1, \SA _2}$ satisfy
\begin{equation}
\gs _{\SA _1 , \SA _2 } \circ (j_1 \Sqcup j_2 ) = \bigl( \RS (j_1) \ot \RS (j_2 )\bigr) \circ \gs  _{\SB _1 , \SB _2} .
\end{equation}
{\bf (b)}
The universal product is associative  iff
\begin{align}
&\bigl( \id _{\RS (\SA _1 ) ^{\ot d}} \ot \RS (\gs _{\SA _2 , \SA _3} )\bigr)  \circ \gs _{\SA _1 , \SA _2 \sqcup \SA _3 }  \nonumber \\
& \ =
\bigl( \RS (\gs _{\SA _1 , \SA _2} )\ot  \id _{\RS (\SA _3 )^{\ot d}}  \bigr)  \circ \gs _{\SA _1 \sqcup \SA _2 , \SA _3 } .
\end{align}
{\bf (c)}
The universal product is unital  iff
\begin{align}
\gs _{\SA _1 , \SA _2 } \circ \iota _1 &= \iota _{\SA _ 1 ^d} ,  \nonumber \\
\gs _{\SA _1 , \SA _2 } \circ \iota _2 &= \iota _{\SA _ 2 ^d} .
\end{align} 
{\bf (d)}  
The universal product is symmetric iff
\begin{equation}
\tau \circ \gs _{\SA _2 ,  \SA _1} =  \gs _{\SA _1 ,  \SA _2}  \circ \tau
\end{equation}
where the first $\tau$ is the tensor switch and the second $\tau$ is the natural identification of $\SA _2 \sqcup \SA _1$ with $\SA _1 \sqcup \SA _2$.
\eprop
\pn
{\it Proof}: Like for $d = m = 1$ in \cite{Lachs1} .
\hfill $\square$
\bn
A \emph{dual group} is a dual semi-group $\SD$ such that there exists an algebra homomorphism
$A$ on $\SD$ with $M \circ (\id \Sqcup A ) \circ \GD = 0 = M \circ (A \Sqcup \id ) \circ \GD $, $M : \SD \sqcup \SD \to \SD$ the multiplication map of the algebra $\SD$.
The free product of two dual groups $(\SD _1 , \GD _1 )$ and $(\SD _2 , \GD _2 )$ is the dual group $\bigl( \SD _1 \sqcup \SD _2 , \tau \circ (\GD _1 \Sqcup \GD _2) \bigr)$ with $\tau$ the identification 
\begin{equation*}
(\SD _1 \sqcup \SD _1 ) \sqcup (\SD _2 \sqcup \SD _2 ) =
(\SD _1 \sqcup \SD _2 ) \sqcup (\SD _1 \sqcup \SD _2 )
\end{equation*}
Similarly, we define the free product of $m$ dual semi-groups and the $m$-fold free product $(\SD ^{\sqcup m} , \GD ^{\sqcup m})$ of a dual semi-group with itself.
For a u.a.u.-product $\odot$ in $\AlgProbLB$ and for $\phi _1 , \phi _2 \in (\SD ^{\sqcup m} )^{*d} = \bigl( (\SD ^{\sqcup m} ) ^ d \bigr)^*$ we define the convolution product of $\phi _1$ and $\phi _2$ as the element of $(\SD ^{\sqcup m} )^{*d}$ given by
\begin{equation}
\label{conv}
\phi _1 \conv \phi _2 := (\phi _1 \odot \phi _2 ) \circ \GD ^{\sqcup m } .
\end{equation}
The space $(\SD ^{\sqcup m} )^{*d}$ is a monoid with respect to convolution.
\par
For the notions of coalgebras, bialgebras and Hopf algebras see, for example, \cite{DNR}.
One shows (cf. \cite{Lachs1})
\bthm
\label{Lachs}
{\bf (a)}
The triplet 
\begin{equation}
\label{bialgebra}
\Bigl( \RS \bigl( (\SD ^{\sqcup m} )^d \bigr) , \RS (\gs _{\SD ^{\sqcup m} , \SD ^{\sqcup m}} \circ \GD ^{\sqcup m} ), \RS (0) \Bigr)
\end{equation}
is a commutative bialgebra and 
\begin{equation}
\label{conviso}
\RS (\phi _1 \conv \phi _2 ) = \RS (\phi _1 ) \ast \RS (\phi _2 )
\end{equation}
where $\ast$ on the left side of (\ref{conviso}) is given by the comultiplication of the coalgebra (\ref{bialgebra}).
\pn
{\bf (b)} If $\SD$ is a dual group with antipode $A$, the bialgebra (\ref{bialgebra}) is a Hopf algebra with antipode $\RS (A)$. 
\ethm
\pn
\bremark
{\rm It has been pointed out by Stephanie Lachs \cite{Lachs1} that Theorem \ref{Lachs} in the case $d = m = 1$ establishes a cotensor functor (= colax monoidal functor) between the tensor category $(\AlgProb , \odot)$ and the tensor category $(\ComUnAlgProb , \otimes )$ where $\ComUnAlgProb$ denotes the category of pairs $(\SA , \phi )$ wit $\SA$ a commutative unital algebra and $\phi$ a normalized linear functional on $\SA$. This result can now readily be extended to the general $d, m$-case to obtain a cotensor functor $\SL _{\odot}$ from $(\AlgProbLB , \odot )$ again to $(\ComUnAlgProb , \otimes )$. Through this functor we can associate with each pair $(\SD , \odot )$ with $\SD$ a u.a.u.-product $\odot$ in $\AlgProbLB$ the commutative bialgebra $\SL _{\odot} (\SD )$ given by (\ref{bialgebra}).
If $\SD$ is a dual group,  $\SL _{\odot} (\SD )$ is a commutative Hopf algebra.
Moreover, the cotensor property of $\SL _{\odot}$ also means  that equation (\ref{conviso}) holds.}
\eremark
\bremark
{\rm
If $\SD$ is co-commutative and $\odot$ symmetric the bialgebra $\SL _{\odot} (\SD )$ is co-commutative.}
\eremark
\bremark
{\rm
Let $\SD$ be a dual group. The convolution (\ref{conv}) defines a group structure on the set $(\SD ^{\sqcup m}) ^{*d}$. This follows from the facts that 
\begin{equation*}
\RS : (\SD ^{\sqcup m}) ^{*d} \to \RS \bigl( (\SD ^{\sqcup m} )^d \bigr) ^*
\end{equation*}
is injective and $\RS (\phi )$ as a homomorphism on a Hopf algebra is invertible; cf. \cite{FrMcK}.
}
\eremark
\pn
A dual semi-group $\SD$ is said to be $\mn$-graded if $\SD$ is an $\mn$-graded algebra and $\GD$ is homogeneous of degree 0. 
\bthm
\label{gradedLachs}
Let $\SD$ be an $\mn$-graded dual semi-group and let $\odot$ be 
a u.a.u.-product on $\AlgProbLB$.
\pn
Then the bialgebra $\SL _{\odot} (\SD )$ of Theorem \ref{Lachs} is an $\mn _0$-graded bialgebra.
\ethm
\pn
{\it Proof}: Follows from Theorem \ref{coeff}.\hfill$\square$
\bremark
{\rm
In fact, the bialgebra $\SL _{\odot} (\SD )$ in Theorem \ref{gradedLachs} is  an $\mn _0$-graded \emph{Hopf} algebra; see Remark \ref{evenHopf}.
}
\eremark
\bexam
{\rm
We consider the case of the dual semi-group 
\begin{equation*}
\bigl( \RT (\SV ), \RT (\iota _1 \oplus \iota _2 ) \bigr)
\end{equation*}
where $\iota _{1/2} : \SV \to \SV \oplus \SV$ denote the embeddings of $\SV$ into the first and the second summand of $\SV \oplus \SV$, i.e. 
\begin{equation*}
\RT (\iota _1 \oplus \iota _2 ) : \RT (\SV ) \to \RT (\SV \oplus \SV ) \cong \RT (\SV ) \sqcup \RT (\SV )
\end{equation*}
is the algebra homomorphism such that
\begin{equation*}
\RT (\iota _1 \oplus \iota _ 2 ) (v) = \iota _1 (v) + \iota _2 (v) ,  \ v \in \SV .
\end{equation*}
}
\eexam
\pn
This actually is a dual \emph{group} with the antipode $\RT (\mu )$, $\mu (v) = - v$, $v \in \SV$.
We also assume the vector space $\SV$ to be $\mn $-graded
\begin{equation*}
\SV = \bigoplus_{n = 1} ^{\infty} \SV ^{(n)} .
\end{equation*}
The tensor algebra $\RT (\SV )$ inherits the grading in the unique way to become an $\mn $-graded dual group.
For example, if $\SV$ is the linear span of $x_1 , \ldots , x_n$ we have $\RT (\SV ) = \mc \lan X_n \ran = \mc \lan x_1 , \ldots , x_n\ran$ and a possible grading  is the usual grade of a non-commutative polynomial in $\mc \lan n\ran$ which comes from the grading with $\SV ^{(1)} = \SV$, $\SV ^{(n)} = \{ 0\}$, $n \neq 1$.
\bthm
\label{Hopf}
Let $\SV$ be an $\mn$-graded vector space. 
\pn
Then
the bialgebra 
$\SL _{\odot} \bigl(\RT (\SV )\bigr)$
is an $\mn _0$-graded Hopf algebra.
\ethm
\pn

\section{{\texorpdfstring{$\mn _0$}{\mn 0}}-graded bialgebras}

We come to some general considerations on $\mn _0$-graded bialgebras.
First let $(\SC , \GD , \gd )$ be an $\mn _0$-graded coalgebra so that we have
\begin{gather*}
\SC = \bigoplus_{n = 0} ^{\infty} \SC ^{(n)} ,  \\
\GD \SC ^{(n)} \subset \bigoplus_{k = 0}^n \bigl( \SC ^{(k)} \ot \SC ^{(n - k)}\bigr).
\end{gather*}
Then the counit $\gd$ vanishes on $\SC ^{(n)}$ for $n \neq 0$. Put
\begin{equation*}
\SC _{\gd} ^* := \{ \ga \in \SC ^*  \mit \ga \rest \SC^{(0)} = \gd \rest \SC ^{(0)} \} .
\end{equation*}
For $c \in \SC ^{(n)} $ and $\ga _1 , \ldots , \ga _{n + 1} \in \SC _{\gd} ^*$ we have
\begin{equation*}
\bigl( (\ga _1 - \gd ) \conv \ldots \conv (\ga _{n + 1} - \gd ) \bigr) (c) = 0 .
\end{equation*}
The subspace
\begin{equation*}
\SC ^{(\leq n )} := \bigoplus_{l = 0}^n \SC ^{(l)}
\end{equation*}
is a sub-coalgebra of $\SC$ and $(\ga _1 - \gd ) \conv \ldots \conv (\ga _{n + 1} - \gd )$ vanishes on $\SC ^{(\leq n )}$.
In particular, $(\ga - \gd ) ^{ \conv (n + 1)}$ vanishes on $\SC ^{(\leq n )}$ for $\ga \in \SC _{\gd} ^*$, i.e. 
\begin{equation*}
(\ga - \gd ) \rest \SC ^{(\leq n )}
\end{equation*}
 is nilpotent. 
Therefore, we may form the logarithmic series
\begin{equation}
\label{ln}
\ln _{\conv} \ga := \sum _{l = 1} ^{\infty} (-1) ^{l + 1} \, \frac{(\ga - \gd )^{\conv l}}{l} ,
\end{equation}
first on $\SC ^{(\leq n )}$, then on the whole of $\SC$ so that (\ref{ln}) can be understood pointwise.
We defined a mapping $\ln _{\conv} : \SC _{\gd} ^* \to \SC ^*$.
\par
For $\gb \in \SC ^*$ we form the convolution exponential
\begin{equation}
\label{convexpo}
\exp _{\conv} {\gb} := \sum_{l = 0}^{\infty} \frac{\gb ^{\conv l}}{l !} \in \SC ^*
\end{equation} 
which converges pointwise by the fundamental theorem on coalgebras; see \cite{DNR} and \cite{SchVo}.
Since $\ln _{\conv} \ga \rest \SC ^{(0)} = 0$, the linear functional $\ln _{\conv} \ga$ must be nilpotent on $\SC ^{(\leq n )}$, and the series (\ref{convexpo}) is a finite sum on $\SC ^{(\leq n )}$ for $\gb = \ln _{\conv} \ga$. It follows from a formal power series argument that we have
\begin{equation*}
\exp _{\conv} ({\ln _{\conv} \ga } ) = \ga
\end{equation*}
for all $\ga \in \SC _{\gd} ^*$. 
As $\exp _{\conv} {\gb } \in \SC _{\gd} ^*$ if $\gb \rest \SC ^{(0)} = 0$, we have
\begin{equation*}
\ln _{\conv} (\exp _{\conv} {\gb }) = \gb
\end{equation*}
for $\gb \in \SC ^*$ with $\gb \rest \SC ^{(0)} = 0$. 
Again by nilpotency and a formal power series argument, we have for $\ga _1 , \ga _2 \in \SC _{\gd} ^*$
\begin{align}
\label{CBH}
&\ln _{\conv} ( \ga _1 \conv \ga _2 ) \nonumber \\
& \  = \RH (\ln _{\conv} \ga _1 , \ln _{\conv} \ga _2 ) \nonumber \\
&  \  = \ln _{\conv} \ga _1 + \ln_{\conv} \ga _2 + \frac12 [\ln _{\conv} \ga _1 , \ln _{\conv} \ga _2 ]_{\conv} \nonumber \\
&\hspace{2cm} + \frac{1}{12} [\ln _{\conv} \ga _1 , [\ln _{\conv} \ga _1 , \ln _{\conv} \ga _2 ]_{\conv} ] _{\conv} + \ldots
\end{align}
where $\RH (x, y ) $ denotes the Campbell-Baker-Hausdorff expansion 
\begin{equation*}
\RH (x, y) = \ln ( \re ^x \, \re ^y ) \in \mc \lan x, y \ran
\end{equation*}
for formal power series.
(Notice that the series  (\ref{CBH}) is again finite on $\SC ^{(\leq n )}$.)
\par
We now assume that $\SB$ is an $\mn _0$-graded Hopf algebra with $\SB ^{(0)} = \mc \ei$. In this case $\ga \in \SB _{\gd} ^*$ means that $\ga$ is a normalized linear functional on $\SC$, i.e. $\ga (\ei ) = 1$.
\bremark 
\label{evenHopf}
{\rm
If $\SB$ is an $\mn _0$-graded bialgebra with 1-dimensional $0$-com\-po\-nent, it automatically is  an ($\mn _0$-graded) Hopf algebra; see e.g. \cite{Man}.
In fact, the antipode $\RA$ is given by
\begin{equation}
\label{antipode1}
\RA = \exp_{\conv} (- \RD ) = \sum_{l = 0} ^{\infty} \frac{(- \RD ) ^{*l}}{l !}
\end{equation}
with
\begin{equation}
\label{antipode2}
\RD = \sum_{l = 1} ^\infty (-1)^{l + 1} \frac{(\id _{\SB} - \gd \ei )^{\conv l}}{l}
\end{equation}
where the series (\ref{antipode1}) and (\ref{antipode2})   are finite sums in $\SB$ if evaluated at a point $b \in \SB$. For a unital algebra homomorphism $\ga$ from $\SB$ to the complex numbers we have $\ga \conv (\ga \circ A) = \gd$  so that $\ga$ is always invertible with respect to convolution.
}
\eremark
A linear functional $\gb$ on an $\mn _0$-graded Hopf algebra $\SB$ is called a $\gd$-derivation if
\begin{equation}
\label{derivation}
\gb (ab) = \gb (a)   \gd (b) + \gd (a) \gb (b)
\end{equation}
for all $a, b \in \SB$.
For $\gb _1 , \gb _2 \in \SB^*$ we put
\begin{equation}
\label{LieBracket}
[\beta _1 , \beta _2 ]_{\conv}  = \beta _1 \conv \beta _2 - \beta _2 \conv \beta _1 .
\end{equation}
\bprop
\label{propderiv}
Let $\SB$ be an $\mn _0$-graded Hopf algebra with $\SB ^{(0)} = \mc \ei$.
Moreover, let  $\ga : \SB \to \mc$ be a unital algebra homomorphism and let
$\beta, \beta _1 , \beta _2 : \SB \to \mc$ be $\gd$-derivations.
\pn
Then
we have
\pn
{\bf (a)} $\ln _{\conv} \ga$ is a $\gd$-derivation.
\pn
{\bf (b)} $\exp_{\conv} \beta$ is a unital algebra homomorphism.
\pn
{\bf (c)} $[\beta _1 , \beta _2 ]_{\conv}$
is a $\gd$-derivation.
\eprop
\pn
{\it Proof}:
{(a)}:
Let $a, b \in \SB$ be homogeneous elements with $\deg (a) + \deg (b) = n$.
Denote by $T$ the restriction of
\begin{equation*}
\bigl( \id \ot (\ga - \gd ) \bigr)  \circ  \GD = \bigl( (\id \ot \ga ) \circ \GD \bigr) - \id
\end{equation*}
to the sub-coalgebra $\SB ^{(\leq n )}$ of $\SB$, which $T$ leaves invariant.
Then
\begin{equation*}
T(ab) = (Ta ) (Tb ) + (Ta) (b) + a T(b)
\end{equation*}
and, more generally,
\begin{equation*}
T ^n (ab) = \bigl( M \circ (T \ot \id + \id \ot T + T \ot T ) ^n \bigr) (a \ot b)
\end{equation*}
with $M : \SB \ot \SB \to \SB$ the multiplication of $\SB$.
Since $T$ is nilpotent we can define
\begin{equation*}
\ln (\id + T) := \sum_{l = 1} ^{\infty} (-1) ^{l + 1} \, \frac{T ^l}{l}
\end{equation*}
and we have
\begin{align*}
\ln (\id + T ) (ab) &= \bigl( M \circ \ln (\id \ot \id + T \ot \id + \id \ot T + T \ot T ) \bigr) (a \ot b)   \\
&= \bigl( M \circ \bigl( \ln (\id + T ) \ot \id + \id \ot \ln (\id + T ) \bigr) \bigr) (a \ot b )  \\
& = \ln (\id + T ) (a) (b) + a  \ln (\id + T )(b)
\end{align*}
which is (\ref{derivation}) if we apply $\gd$. 
\pn
{(b)}:
That $\exp _{\conv} \beta$ is a unital algebra homomorphisms  can be seen in a similar way. 
\pn
{(c)}: This follows from the fact that the commutator of two derivations is again a derivation.\hfill$\square$
\bthm
\label{Lie}
Let $\SB$ be
an $\mn _0$-graded Hopf algebra with $\SB ^{(0)} = \mc \ei$.
\pn
Then we have
\pn
{\bf (a)}
The set 
\begin{equation*}
\RG (\SB ) = \{ \ga \in \SB ^* \mit \ga \ \mbox{is a unital algebra homomorphism} \}
\end{equation*}
forms a group with respect to convolution.
\pn
{\bf (b)}
The set 
\begin{equation*}
\RL (\SB ) = \{ \gb \in \SB ^* \mit \gb \ \mbox{is a} \  \gd\mbox{-derivation} \}
\end{equation*}
forms a Lie algebra with respect to the Lie bracket given by (\ref{LieBracket}).
\ethm
\bremark
{\rm
The mappings $\exp_{\conv} : \RL (\SB ) \to \RG (\SB )$ and $\ln_{\conv} : \RG (\SB ) \to \RL (\SB )$ are inverses of each other.
$\exp _{\conv}$ plays the role of the exponential mapping of a Lie group; cf. \cite{FrMcK}.
}
\eremark

\section{Cumulants}
\label{cumulants}
Now we apply the above considerations on $\mn _0$-graded Hopf algebras to the case of the $\mn _0$-graded Hopf algebra $\SL _{\odot} \bigl(\RT (\SV )\bigr)$ of Theorem \ref{Hopf} 
coming from a u.a.u.-product $\odot$ in the category $\AlgProbLB$ and a vector space $\SV$, which is considered as $\mn$-graded with $\SV ^{(1)} = \SV$.
Put $\SE = \bigl( \RT (\SV ) ^{\sqcup m} \bigr) ^ d$.
For $\phi \in \SE ^* = \bigl( \RT (\SV ) ^{\sqcup m } \bigr) ^{* d}$ we have that $\RS (\phi ) \in \RS (\SE ) ^*$ is an algebra homomorphism. 
Therefore, by Proposition \ref{propderiv}, $\ln _{\conv} \bigl( \RS (\phi ) \bigr)$ is an $\RS (0)$-derivation on $\RS (\SE )$. 
Such a derivation $D$ is of the from $D = D (\psi )$ with $\psi = D \rest \SE$. Here 
we put
\begin{equation*}
D (\psi ) \rest \SE  ^{\ot _s n} =  \left\{ \begin{array}{ll} \psi & \mbox{if} \ n = 1  \\
0 & \mbox{otherwise}
\end{array} \right.
\end{equation*}
for $\psi \in \SE ^*$.
We call
\begin{equation*}
\ln_{\odot} (\phi ) := \ln_{\conv} \bigl( \RS (\phi )\bigr) \rest \SE \in \SE ^*
\end{equation*}
the cumulant functional of $\phi$ and we call
\begin{equation*}
\ln _{\odot} : \SE ^* \to \SE ^*
\end{equation*}
the cumulant mapping of $\odot$ for the vector space $\SV$.
\bthm
Let $\phi _1 , \phi _2 \in \SE ^*$.
\pn
Then we have
\begin{equation}
\label{CBH2}
\ln _{\odot} (\phi _1 \conv \phi _2 )  = \RH \bigl( D(\ln _{\odot} (\phi _1) ), D(\ln _{\odot} (\phi _2 ) )\bigr) \rest \SE .
\end{equation}
If the underlying universal product $\odot$ is symmetric, we have
\begin{equation}
\label{additivity}
\ln _{\odot} (\phi _1 \conv \phi _2 ) = \ln _{\odot} (\phi _1 ) + \ln _{\odot} (\phi _2 )
\end{equation}
\ethm
\pn
{\it Proof}:
The first part is (\ref{CBH}).
The second part follows from $D(\psi _1 ) + D(\psi _2 ) = D(\psi _1 + \psi _2 )$ and the fact that convolution in $\RS (\SE ) ^*$ is commutative if $\odot$ is symmetric.\hfill$\square$
\bexam
{\rm The products given by the tensor, free, Boolean, c-free and bi-free independence are symmetric, so that in these cases we have the well-known additivity (\ref{additivity}) of cumulants; cf. for example \cite{FrMcK}. The monotonic independence of Muraki is not symmetric, which means that there may be terms of order $\geq 2$ in the Campbell-Baker-Hausdorff expansion (\ref{CBH2}); see again \cite{FrMcK}.
}
\eexam
\bremark
{\rm
As we have seen, $\mc$-valued  unital algebra homomorphisms and $\S (0)$-derivations of $\SL _{\odot} \bigl(\RT (\SV )\bigr)^*$ both are given by an element of $\SE ^*$.
We will write $\SE _{\SV} ^*$ for $\SE ^*$ to emphasize the dependence on $\SV$.
Therefore, $\RG \bigl( \SL _{\odot} \bigl(\RT (\SV )\bigr)\bigr)$ turns  $\SE _{\SV } ^*$ into a group and $\RL \bigl( \SL _{\odot} \bigl(\RT (\SV )\bigr)\bigr)$ turns  $\SE _{\SV} ^*$ into a Lie algebra. In this sense, the cumulant mapping $\ln _{\odot}$ is the restriction of $\ln _{\conv}$ to $\SE _{\SV} ^*$ and the restriction of $\exp _{\conv}$ to $\SE _{\SV} ^*$ is the inverse of $\ln _{\odot}$.
Let us denote this inverse mapping by $\exp _{\odot}$.
}
\eremark
\par
We have seen (Theorems \ref{Hopf} and \ref{Lie}) that a u.a.u.-product on $\AlgProbLB$ gives rise to Lie algebras $\RL \bigl( \SL _{\odot}(\RT (\SV )) \bigr) = \SE _{\SV} ^*$ where $\SV$ runs through all vector spaces. 
The following considerations show that the product $\odot$ can be reconstructed from this family of Lie algebras if we also know the cumulant mappings,
which means that $\odot$ is determined by its cumulants and cumulant Lie algebras.
First, the convolution of $\phi _1 , \phi _2 \in \SE _{\SV} ^*$ is given by
\begin{equation}
\label{equation1}
\phi _1 \conv \phi _2  = \exp _{\odot} \bigl( \RH (\ln _{\odot } \phi _1, \ln _{\odot} \phi _2 ) \bigr) .
\end{equation}
The right hand side of (\ref{equation1}) can be calculated from the cumulant mappings and cumulant Lie algebras, i.e. we have
the convolutions (\ref{conv}) of $d$-tuples of linear functionals on the dual groups $\RT (\SV )^{\sqcup m} = \RT (\SV ^m )$ associated with $\odot$.
Now let 
\begin{equation*}
\bigl(\SA _1 , (\SA _1 ^{(1)} , \ldots , \SA _1 ^{(m)} ), \phi _1 \bigr), \bigl(\SA _2 , (\SA _2 ^{(1)} , \ldots , \SA _2 ^{(m)} ), \phi _2 \bigr)
\end{equation*}
 be two objects in $\AlgProbLB$.
Put
\begin{equation*}
j_i ^{(l)} = \RT (\id _{\SA ^{(l)} _i }) : \RT (\SA _i ^{(l)} ) \to \SA _ i ^{(l)}
\end{equation*} 
$i = 1, 2$, $l \in [m]$, equal to 
the algebra homomorphism given by the multiplication in $\SA _i ^{(l)}$.
Then
\begin{align}
\label{equation2}
\phi _1 \odot \phi _2    &= 
( \phi _1  \odot  \phi _2 ) \circ    \bigl( (j _1 \sqcup 0 ) \Sqcup (0 \sqcup j _2) \bigr) \circ \GD ^{\Sqcup m}\rest (\SA _1 \sqcup \SA _2 ) ^d 
\\
\nonumber
&=  \bigl( \phi _1 \circ (j _1 \sqcup 0 )\bigr) \odot \bigl( \phi _2 \circ (0 \sqcup j _2 )\bigr) \circ \GD ^{\Sqcup m}\rest (\SA _1 \sqcup \SA _2 ) ^d 
\\
\nonumber
&= \bigl( \phi _1 \circ (j _1 \sqcup 0 )\bigr) \conv \bigl( \phi _2 \circ (0 \sqcup j _2 )\bigr) \rest (\SA _1 \sqcup \SA _2 ) ^d 
\subset \bigl( \RT (\SV ^m )\bigr) ^d 
\end{align}
with $\conv$ the convolution of $d$-tuples of linear functionals on the dual group 
$\RT (\SV ^m )$, $\SV = \SA _1 ^{(1)} \oplus \SA _2 ^{(1)} \oplus \ldots \oplus
\SA _1 ^{(m)} \oplus \SA _2 ^{(m)}$.
Since an element of
\begin{equation*}
\SA _1 \sqcup \SA _2 = (\SA _1 ^{(1)} \sqcup \ldots \sqcup \SA _1 ^{(m)} ) \sqcup (\SA _2 ^{(1)} \sqcup \ldots \sqcup \SA _2 ^{(m)} )
\end{equation*}
lies in a sub-algebra generated by a \emph{finite} number of elements of $\SA _i ^{(l)}$, $i = 1, 2$, $l \in [m]$, a u.a.u.-product can be reconstructed from the Lie algebras and the cumulants in the cases $\SV = \mc ^n$, $n \in \mn$, that is $\odot$ can be reconstructed from the \lq $n$-th order cumulants\rq \ and the \lq $n$-th order cumulant Lie algebras\rq \
where $n$ runs through all natural numbers.
Notice that
\begin{align*}
\RL \bigl( \SL _{\odot} \bigl(\RT (\mc ^n) )\bigr)\bigr)  
&= \bigl( \RT (\mc ^n ) ^{\sqcup m} \bigr) ^{*d}  \\
&= \bigl( \mc \langle X_n \rangle ^{\sqcup m} \bigr) ^{*d}  \\
&= \bigl( \mc \langle X_{nm} \rangle  \bigr) ^{*d}  \\
&= \mc \langle \langle X_{nm} \rangle \rangle ^d .
\end{align*}
Put $\RL _{\odot , n} = \mc \langle \langle X_{nm} \rangle \rangle ^d$ and denote by
\begin{equation*}
\ln _{\odot , n } : \mc \langle \langle X_{nm} \rangle \rangle ^d  \to \mc \langle \langle X_{nm} \rangle \rangle ^d 
\end{equation*}
the cumulant mapping for $\mc ^n$. 
 We summarize
\bthm
A
u.a.u.-product $\odot$ can be reconstructed from the series
\begin{equation*}
(\ln _{\odot, n } , \RL _{\odot , n} )_{n \in \mn}.
\end{equation*}
\pn
More precisely, $\odot$ is given by (\ref{equation2}) with the convolution $\conv$ on 
$\mc \langle \langle X_{nm} \rangle \rangle ^{*d}$ defined by (\ref{equation1}).
\ethm
\bn
{\bf Acknowledgement:}
The authors thank Stephanie Lachs and Malte Gerhold for their  interest and valuable comments and contributions.

\bibliographystyle{alpha}
\bibliography{mybib-2}

\newcommand{\Swap}[2]{#2#1}\newcommand{\Sort}[1]{}
\begin{thebibliography}{AHLV15}

\bibitem[AHLV15]{AHLV}
Octavio Arizmendi, Takahiro Hasebe, Franz Lehner, and Carlos Vargas.
\newblock Relations between cumulants in noncommutative probability.
\newblock {\em Infin. Dimens. Anal. Quantum Probab. Relat. Top.}, 282:56--92,
  2015.

\bibitem[AM10]{AgMa}
Marcelo Aguiar and Swapneel Mahajan.
\newblock {\em Monoidal functors, species and {H}opf algebras}, volume~29 of
  {\em CRM Monograph Series}.
\newblock American Mathematical Society, Providence, RI, 2010.

\bibitem[Ans11]{Ansh}
Michael Anshelevich.
\newblock Two-state free {B}rownian motions.
\newblock {\em J. Funct. Anal.}, 260(2):541--565, 2011.

\bibitem[BGS02]{BGS02}
Anis Ben~Ghorbal and Michael Sch{\"u}rmann.
\newblock {Non-com\-mu\-ta\-tive no\-tions of sto\-chas\-tic
  in\-de\-pen\-dence}.
\newblock {\em Math. Proc. Cam\-bridge Philos. Soc.}, 133(3):531--561, 2002.

\bibitem[BGS05]{BGS05}
Anis Ben~Ghorbal and Michael Sch{\"u}rmann.
\newblock Quantum {L}\'evy processes on dual groups.
\newblock {\em Math. Z.}, 251(1):147--165, 2005.

\bibitem[BLS96]{BLSp}
Marek Bo{\.z}ejko, Michael Leinert, and Roland Speicher.
\newblock Convolution and limit theorems for conditionally free random
  variables.
\newblock {\em Pacific J. Math.}, 175(2):357--388, 1996.

\bibitem[BS91]{BSp}
Marek Bo{\.z}ejko and Roland Speicher.
\newblock {$\psi$}-independent and symmetrized white noises.
\newblock In {\em Quantum probability \& related topics}, QP-PQ, VI, pages
  219--236. World Sci. Publ., River Edge, NJ, 1991.

\bibitem[DNR01]{DNR}
Sorin D{\u{a}}sc{\u{a}}lescu, Constantin N{\u{a}}st{\u{a}}sescu, and
  {\c{S}}erban Raianu.
\newblock {\em Hopf algebras}, volume 235 of {\em Monographs and Textbooks in
  Pure and Applied Mathematics}.
\newblock Marcel Dekker, Inc., New York, 2001.

\bibitem[FM15]{FrMcK}
Roland Friedrich and John McKay.
\newblock Homogeneous {L}ie groups and quantum probability.
\newblock {\em arXiv:1506.07089v1}, 2015.

\bibitem[Fra06]{Fra06}
Uwe Franz.
\newblock L\'evy processes on quantum groups and dual groups.
\newblock In {\em Quantum independent increment processes {II}}, volume 1866 of
  {\em Lecture Notes in Math.}, pages 161--257. Springer, Berlin, 2006.

\bibitem[Ger14]{Ger1}
Malte Gerhold.
\newblock {\em On several problems in the theory of comonoidal systems and
  subproduct systems}.
\newblock PhD thesis, Greifswald, 2014.

\bibitem[GL15]{GeLa}
Malte Gerhold and Stephanie Lachs.
\newblock Classification and {GNS}-construction for general universal products.
\newblock {\em Infin. Dimens. Anal. Quantum Probab. Relat. Top.},
  18(1):1550004, 29, 2015.

\bibitem[Has11]{Has}
Takahiro Hasebe.
\newblock Conditionally monotone independence {I}: {I}ndependence, additive
  convolutions and related convolutions.
\newblock {\em Infin. Dimens. Anal. Quantum Probab. Relat. Top.},
  14(3):465--516, 2011.

\bibitem[HP84]{HudPar}
R.~L. Hudson and K.~R. Parthasarathy.
\newblock Quantum {I}to's formula and stochastic evolutions.
\newblock {\em Comm. Math. Phys.}, 93(3):301--323, 1984.

\bibitem[HS11]{Has2}
Takahiro Hasebe and Hayato Saigo.
\newblock The monotone cumulants.
\newblock {\em Ann. Inst. Henri {Poincar\'e} Probab. Stat.}, 47(4):1160--1170,
  2011.

\bibitem[Lac15]{Lachs1}
Stephanie Lachs.
\newblock {\em A new family of universal products and aspects of a non-positive
  quantum probability theory}.
\newblock PhD thesis, Greifswald, 2015.

\bibitem[Lu97]{Lu}
Y.~G. Lu.
\newblock An interacting free {F}ock space and the arcsine law.
\newblock {\em Probab. Math. Statist.}, 17(1, Acta Univ. Wratislav. No.
  1928):149--166, 1997.

\bibitem[Man08]{Man}
Dominique Manchon.
\newblock Hopf algebras in renormalisation.
\newblock In {\em Handbook of algebra. {V}ol. 5}, volume~5 of {\em Handb.
  Algebr.}, pages 365--427. Elsevier/North-Holland, Amsterdam, 2008.

\bibitem[ML98]{McL}
Saunders Mc~Lane.
\newblock {\em Categories for the working mathematician}, volume~5 of {\em
  Graduate Texts in Mathematics}.
\newblock Springer-Verlag, New York, second edition, 1998.

\bibitem[Mur97]{Mur97}
Naofumi Muraki.
\newblock {Noncommutative Brownian motion in monotone Fock space}.
\newblock {\em Commun.\ Math.\ Phys.}, 183:557--570, 1997.

\bibitem[Mur01]{Mur01}
Naofumi Muraki.
\newblock {Monotonic independence, monotonic central limit theorem and
  monotonic law of small numbers}.
\newblock {\em Infin.\ Dimens.\ Anal.\ Quantum Probab.\ Relat.\ Top.},
  4:39--58, 2001.

\bibitem[Mur02]{Mur02}
Naofumi Muraki.
\newblock The five independences as quasi-universal products.
\newblock {\em Infin. Dimens. Anal. Quantum Probab. Relat. Top.},
  5(1):113--134, 2002.

\bibitem[Mur03]{Mur03}
Naofumi Muraki.
\newblock The five independences as natural products.
\newblock {\em Infin. Dimens. Anal. Quantum Probab. Relat. Top.},
  6(3):337--371, 2003.

\bibitem[Mur13]{Mur4}
Naofumi Muraki.
\newblock A simple proof of the classification theorem for positive natural
  products.
\newblock {\em Probab. Math. Statist.}, 33(2):315--326, 2013.

\bibitem[Spe90]{Sp}
Roland Speicher.
\newblock A new example of ``independence'' and ``white noise''.
\newblock {\em Probab. Theory Related Fields}, 84(2):141--159, 1990.

\bibitem[Spe97]{Spe97}
Roland Speicher.
\newblock On universal products.
\newblock In {\em Free probability theory ({W}aterloo, {ON}, 1995)}, volume~12
  of {\em Fields Inst. Commun.}, pages 257--266. Amer. Math. Soc., Providence,
  RI, 1997.

\bibitem[SV14]{SchVo}
Michael Sch{\"u}rmann and Stefan Vo{\ss}.
\newblock Schoenberg correspondence on dual groups.
\newblock {\em Comm. Math. Phys.}, 328(2):849--865, 2014.

\bibitem[Voi85]{Voi85}
Dan-Virgil Voiculescu.
\newblock Symmetries of some reduced free product {$C^\ast$}-algebras.
\newblock In {\em Operator algebras and their connections with topology and
  ergodic theory ({B}u\c steni, 1983)}, volume 1132 of {\em Lecture Notes in
  Math.}, pages 556--588. Springer, Berlin, 1985.

\bibitem[Voi87]{Voi87}
Dan-Virgil Voiculescu.
\newblock {Dual algebraic structures on operator algebras related to free
  products}.
\newblock {\em J.\ Operator Theory}, 17:85--98, 1987.

\bibitem[Voi14]{Voi123}
Dan-Virgil Voiculescu.
\newblock Free probability for pairs of faces {I}.
\newblock {\em Comm. Math. Phys.}, 332(3):955--980, 2014.

\bibitem[vW73]{vW1}
Wilhelm von Waldenfels.
\newblock {An approach to the theory of pressure broadening of spectral lines}.
\newblock In M.~Behara, K.~Krickeberg, and J.~Wolfowitz, editors, {\em
  Probability and Information Theory II}, number 296 in Lect.\ Notes Math.
  Springer, 1973.

\end{thebibliography}

\enddocument